%% file: arxiv_manuscript.tex
\documentclass[11pt]{article}

\usepackage{amsfonts}
\usepackage{amssymb}
\usepackage{booktabs}
\usepackage[margin=15pt,font=small,labelfont={bf,sf},justification=justified]{caption}
\usepackage[superscript,biblabel]{cite}
\usepackage{color}
\usepackage{float}
\usepackage{graphicx}
\usepackage{lscape}
\usepackage{mathtools}
\usepackage{multibib}
\usepackage{multirow}
\usepackage{tabularx}
\usepackage{titleref}
\usepackage{url}
\usepackage{subfig}
\usepackage{bm}
\usepackage{IEEEtrantools}
\usepackage{placeins}
\usepackage{etoolbox}
\usepackage{todonotes}

\usepackage[scaled=.98,p]{XCharter}
\usepackage[scaled=1.04,varqu,varl]{inconsolata}
\usepackage[type1]{cabin}
\usepackage[uprightscript,xcharter,vvarbb,scaled=1.05]{newtxmath}
\usepackage[cal=euler,scr=boondoxo,bb=boondox]{mathalpha}  
\usepackage[T1]{fontenc}
\usepackage[final,protrusion=true,expansion=true]{microtype}

\usepackage{titling}
\usepackage{authblk} 
\pretitle{\begin{center}\bfseries\sffamily\LARGE}
\posttitle{\end{center}}
\usepackage{abstract}

\usepackage{sectsty} 
\allsectionsfont{\sffamily}

\usepackage[left=1in,right=1in,top=1in,bottom=1in,includefoot,heightrounded]{geometry}
\usepackage[hidelinks]{hyperref}

\newcites{supp}{SUPPLEMENTAL REFERENCES}

\makeatletter
\patchcmd{\LS@rot}{90}{-90}{}{}
\patchcmd{\endlandscape}{90}{-90}{}{}
\makeatother

\newcommand{\beginsupplement}{%
        \setcounter{table}{0}
        \renewcommand{\thetable}{S\arabic{table}}%
        \renewcommand{\theHtable}{supp.\arabic{table}}%
        \setcounter{figure}{0}
        \renewcommand{\thefigure}{S\arabic{figure}}%
        \renewcommand{\theHfigure}{supp.\arabic{figure}}%
        \setcounter{section}{0}
        \renewcommand{\thesection}{S\arabic{section}}
        \renewcommand{\theHsection}{supp.\arabic{section}}%
        \setcounter{equation}{0}
        \renewcommand{\theHequation}{supp.\arabic{equation}}%
        \renewcommand{\theequation}{S\arabic{equation}}     
}

\input{shared_macros}
\ieeeformatfalse
\disableequationbreaks

\newcommand{\resultfigwidth}{0.35\textwidth}
\newcommand{\resultfigcolsep}{\hspace{1em}}
\newcommand{\resultfigrowsep}{\\}
\newcommand{\resultfloatbarrier}{\FloatBarrier}

\title{Composite B-Spline Current Deposition and Interpolation Operators for Thin-Wire Finite-Difference Time-Domain Simulations}
\author[1]{Cole Gruninger}
\author[2]{Boyce E. Griffith}
\affil[1]{Department of Mathematics, Fordham University}
\affil[2]{Department of Mathematics, University of North Carolina at Chapel Hill}
\date{}

\begin{document}
\maketitle

\begin{abstract}
\noindent\input{manuscript_abstract}
\end{abstract}


\input{manuscript_body}

\renewcommand\refname{REFERENCES}
\bibliographystyle{unsrturl}
\bibliography{references}

\clearpage
\beginsupplement
\section{Discrete Energy Conservation of the Leapfrog FDTD--Thin-Wire System}
\label{app:energy}
\input{supplementary_material_body}

\end{document}

%% file: shared_macros.tex
\newcommand{\vE}{\mathbf{E}}
\newcommand{\vF}{\mathbf{F}}
\newcommand{\vG}{\mathbf{G}}
\newcommand{\vH}{\mathbf{H}}
\newcommand{\vJ}{\mathbf{J}}
\newcommand{\vx}{\mathbf{x}}
\newcommand{\vX}{\mathbf{X}}
\newcommand{\energy}{\mathcal{E}}
\newcommand{\vtau}{\boldsymbol{\tau}}
\newcommand{\ds}{\Delta s}
\newcommand{\dx}{\Delta x}
\newcommand{\dy}{\Delta y}
\newcommand{\dz}{\Delta z}
\newcommand{\dt}{\Delta t}
\newcommand{\Scal}{\mathcal{S}}
\newcommand{\Ical}{\mathcal{I}}
\newcommand{\R}{\mathbb{R}}
\newcommand{\divh}{\nabla_h \cdot}
\newcommand{\curlh}{\nabla_h \times}

\newif\ifnarrowequations
\newcommand{\enableequationbreaks}{%
  \narrowequationstrue
  \def\eqbreak{\nonumber\\}%
  \def\eqcontinue{&}%
  \def\eqcontinuequad{&\quad}%
  \def\eqfirstlhs##1{&##1}%
  \def\eqfirststep{\nonumber\\}%
  \def\eqfirstrel{={}&}%
  \def\eqrel{={}&}%
  \def\eqforcedbreak{\nonumber\\}%
  \def\eqforcedcontinue{&}%
  \def\eqieeebreak{\nonumber\\}%
  \def\eqieeecontinuequad{&& \quad\,}%
}
\newcommand{\disableequationbreaks}{%
  \narrowequationsfalse
  \def\eqbreak{\relax}%
  \def\eqcontinue{\relax}%
  \def\eqcontinuequad{\relax}%
  \def\eqfirstlhs##1{##1}%
  \def\eqfirststep{\relax}%
  \def\eqfirstrel{&=}%
  \def\eqrel{&=}%
  \def\eqforcedbreak{\nonumber\\}%
  \def\eqforcedcontinue{&\quad}%
  \def\eqieeebreak{\relax}%
  \def\eqieeecontinuequad{\relax}%
}
\enableequationbreaks

\newif\ifieeeformat
\ieeeformatfalse
\newcommand{\formatparstart}[2]{\ifieeeformat\IEEEPARstart{#1}{#2}\else #1#2\fi}

%% file: manuscript_abstract.tex
Holland--Simpson thin-wire finite-difference time-domain (FDTD) simulations of obliquely oriented closed-loop antennas suffer from persistent low-frequency parasitic currents caused by the failure of the current-deposition operator to conserve charge. This deposition operator, together with an interpolation operator that samples the tangential component of the electric field along the wire, can be realized as regularizations of distributions: the wire current is deposited as a source term by integrating it against a regularized delta function along the wire, and the electric field is sampled back to the wire by integrating it against the same regularized delta function supported on the wire. We show that charge conservation requires the deposited current to be discretely divergence-free if the wire carries a constant current, and interpret this condition in a distributional context. We then introduce a family of composite B-spline regularizations $\delta_h^{(n)}$ that satisfy this charge-conservation condition to machine precision, provided that the coupling line integrals are computed exactly. Such exact evaluation is achievable in this setting because the B-spline kernels are piecewise polynomial and their breakpoints along the wire are known a priori, so composite Gauss--Legendre quadrature with subinterval breakpoints placed at every grid-plane crossing computes the coupling integrals to machine precision. Taking the interpolation operator to be the discrete adjoint of the current-deposition operator keeps the spatial discretization skew-symmetric and ensures that a discretely irrotational electric field drives no net electromotive force around a closed wire loop. Numerical experiments on a center-fed dipole and circular and square loop antennas at three orientations relative to the grid demonstrate that the charge-conserving composite B-spline regularized delta functions produce orientation-independent impedance values in agreement with the known impedance characteristics of these antennas. By contrast, the standard trilinear regularization produces parasitic low-frequency currents in the closed-loop antennas, rendering the extracted impedance unphysical.

%% file: manuscript_body.tex
\section{Introduction}
\label{sec:intro}

\formatparstart{T}{he} Holland and Simpson thin-wire formalism~\cite{holland1981} permits wires of subgrid radius to be embedded in finite-difference time-domain (FDTD) simulations on the Yee grid~\cite{yee1966,taflove2005} without resolving the wire cross-section on the three-dimensional mesh. The wire is modeled as a one-dimensional transmission line, and its current and voltage evolve according to the telegrapher equations. The coupling between the wire and the surrounding electromagnetic field is realized through two operators: a \emph{current deposition} operator that deposits wire current onto the Yee grid as a source current density, and an \emph{interpolation} operator that samples the tangential electric field from the grid back to the wire. Both operators can be realized as regularizations of distributions~\cite{schwartz1998}. The current deposition operator represents the wire current as a distributional source term, and the interpolation operator pairs the electric field with the same distribution to extract its tangential component along the wire.

A well-known difficulty arises in thin-wire FDTD simulations of loop antennas if the wire is oriented obliquely to the Cartesian grid axes: the simulated current can develop persistent low-frequency spurious oscillations that are absent from the physical solution~\cite{edelvik2003,berenger2019,guiffaut2010,ledfelt2002}. In the lossless, source-free setting, the leapfrog time integrator used in standard FDTD is symplectic and conserves a discrete analog of the electromagnetic energy~\cite{schuhmann2001,taflove2005}, making the method nondissipative and free of artificial damping. This property is consistent with energy conservation in Maxwell's equations and is beneficial for long-time wave propagation, but it also means that any parasitic modes introduced by the spatial discretization are not damped by the time integrator and persist indefinitely.

An early work by Ledfelt~\cite{ledfelt2002} identified the orientation dependence and persistent low-frequency oscillations produced by a simple trilinear regularization of the delta function. That work introduced a shell-averaged interpolation scheme in which the electric field is averaged over a small circular shell of points around the wire axis, and distributed the current back to the grid through the discrete adjoint of this interpolation. The resulting spatial discretization is skew-symmetric, which, combined with the leapfrog time integrator, yields discrete energy conservation. Closed loops nonetheless still exhibited a weak undamped low-frequency oscillation, which that work attributed to overlap between the support shells of adjacent wire segments. Subsequent work by Edelvik~\cite{edelvik2003} refined this approach using a radially symmetric cosine kernel supported on a cylindrical tube around the wire axis, again taking the current deposition operator as the discrete adjoint of the interpolation operator. Neither work identified the underlying source of the parasitic oscillations in the current.

The mechanism behind these residual oscillations was later characterized in the work of Guiffaut, Reineix, and Pecqueux~\cite{guiffaut2010,guiffaut2012}, which showed that the current deposition operator deposits a nonphysical charge density on the grid that drives the parasitic currents. That work addressed the problem with a regularized delta function that is piecewise constant along the wire's dominant Cartesian direction and piecewise linear in the two transverse directions, combined with a geometric trace-continuity rule at each cell node. Subsequent work by B\'{e}renger~\cite{berenger2019} gave a more systematic analysis, tracing the parasitic solutions to an imbalance of current entering and exiting fictitious \emph{Q-cells} surrounding the nodes of the Yee grid. That work derived a charge-conservation criterion from the case of constant wire current and presented a geometric implementation satisfying it, but did not pair the current deposition operator with an adjoint interpolation. While the scheme remained stable under the standard FDTD CFL condition, the lack of skew-symmetry permits spurious energy creation or destruction. A more recent review by Qiao et al.~\cite{qiao2023} tested several charge-conservation schemes adapted from the particle-in-cell literature, including the Villasenor--Buneman volume-weighting scheme~\cite{villasenor1992} and its Zigzag and Esirkepov variants.

In this paper we reinterpret B\'{e}renger's condition in a distributional context and propose an alternative implementation based on composite B-spline regularizations of the delta function~\cite{gruninger2026,li2025}. The central observation is that a closed wire loop carrying a constant current produces a source term that is, in a distributional sense, a curl. Charge conservation in the thin-wire FDTD system then requires that this distributional curl be mapped to a discrete curl on the Yee grid. This reinterpretation unifies B\'{e}renger's geometric criterion and the Guiffaut--Reineix--Pecqueux trace-continuity condition within a single framework, and yields current deposition and interpolation operators that are skew-symmetric by construction, restoring the reciprocity that B\'{e}renger's original scheme lacked.

This distributional perspective leads to three primary contributions. We identify the discrete analog of the distributional divergence-free property of a constant wire current, expressed in condition~\eqref{eq:berenger_cond}, and propose a family $\delta_h^{(n)}$ of regularized delta functions indexed by an integer $n \ge 0$ that satisfies this condition and yields a systematic hierarchy of charge-conserving choices of increasing smoothness. Charge conservation in this construction holds to machine precision as long as the line integrals associated with the current-deposition operator are computed exactly, which we achieve through composite Gauss--Legendre quadrature with breakpoints at the grid-plane crossings induced by the chosen kernel pair. The interpolation operator is then constructed as the discrete adjoint of the current-deposition operator, so that a discretely irrotational electric field drives no electromotive force around a closed wire loop, just as in the continuous setting.

\section{Background}
\label{sec:background}

\subsection{The Yee Grid and the FDTD Method}

The FDTD method discretizes Maxwell's curl equations on a staggered Yee grid~\cite{yee1966,taflove2005}. On a uniform grid with spacings $\dx$, $\dy$, $\dz$, the electric-field components are located at edge midpoints:
\begin{align}
E_x &\text{ at } \bigl((i\!+\!\tfrac{1}{2})\dx,\; j\dy,\; k\dz\bigr), \\
E_y &\text{ at } \bigl(i\dx,\; (j\!+\!\tfrac{1}{2})\dy,\; k\dz\bigr), \\
E_z &\text{ at } \bigl(i\dx,\; j\dy,\; (k\!+\!\tfrac{1}{2})\dz\bigr),
\end{align}
while the magnetic-field components are located at face centers with complementary staggering. The leapfrog updates for the magnetic and electric fields are
\begin{align}
\vH^{n+1/2} &= \vH^{n-1/2} - \frac{\dt}{\mu_0} \curlh \vE^n, \label{eq:Hupdate}\\
\vE^{n+1} &= \vE^n + \frac{\dt}{\epsilon_0}
  \bigl(\curlh \vH^{n+1/2} - \vJ^{n+1/2}\bigr), \label{eq:Eupdate}
\end{align}
in which $\curlh$ is the discrete curl on the Yee grid and $\vJ^{n+1/2}$ is the source current density deposited onto the grid by the thin-wire current deposition operator, defined in Section~\ref{sec:operators}.

A fundamental identity of the Yee finite difference stencils is that the discrete divergence of a discrete curl vanishes identically:
\begin{equation}
\divh (\curlh \mathbf{A}) = 0 \quad \text{for any grid field } \mathbf{A}.
\label{eq:div_curl_zero}
\end{equation}
This identity is not merely a coincidence: it reflects the fact that the Yee grid discretizes Maxwell's equations in a manner consistent with their natural differential form interpretation, in which the electric field is a 1-form and the magnetic field a 2-form. Yee's original work did not acknowledge this structure explicitly, and it appears to have been a fortunate consequence of the staggered grid construction. Interestingly, Lebedev~\cite{lebedev1964,lebedev1964a} had independently arrived at finite difference discretizations preserving this same orthogonal structure just a few years prior to Yee's publication. These properties have since been systematically understood through the frameworks of mimetic finite differences~\cite{hyman1999,hyman1999a} and finite element exterior calculus (FEEC)~\cite{arnold2010,arnold2018,glasser2023,bossavit1999,bossavit2000}, in which the Yee grid corresponds to the lowest-order Whitney form discretization of the discrete version of the de Rham complex. A direct consequence is the discrete Helmholtz decomposition: any discrete vector field on the Yee grid can be decomposed as
\begin{equation}
\vE_h = \nabla_h \phi + \curlh \mathbf{W} + \boldsymbol{\gamma},
\label{eq:helmholtz}
\end{equation}
in which $\phi$ is a discrete scalar potential, $\mathbf{W}$ is a discrete vector potential, and $\boldsymbol{\gamma}$ is a discrete harmonic component. The three subspaces are mutually orthogonal with respect to the natural $\ell^2$ inner product on the Yee grid, and in particular $\langle \nabla_h \phi,\, \curlh \mathbf{W} \rangle_h = 0$ for any $\phi$ and $\mathbf{W}$. Equation~\eqref{eq:div_curl_zero} is precisely the algebraic statement underlying this orthogonality, and the two facts are equivalent.

The leapfrog scheme~\eqref{eq:Hupdate}--\eqref{eq:Eupdate} is a symplectic integrator that preserves a discrete analog of Poynting's theorem~\cite{taflove2005,schuhmann2001}. Applying $\divh$ to~\eqref{eq:Eupdate} and using~\eqref{eq:div_curl_zero} gives
\begin{equation}
\epsilon_0 \divh \vE^{n+1}
= \epsilon_0 \divh \vE^n - \dt\, \divh \vJ^{n+1/2}.
\label{eq:gauss_update}
\end{equation}
Gauss's law $\epsilon_0 \divh \vE = \rho_h$ is therefore preserved exactly if $\divh \vJ = 0$, a consequence of the spatial identity~\eqref{eq:div_curl_zero} alone. The energy-conserving nature of the leapfrog scheme, however, means that any spurious charge introduced by a non-divergence-free source current is not damped by the time integrator --- it persists indefinitely and drives a growing contamination of the electric field. Charge conservation in the source current is therefore not merely a consistency requirement; it is essential for the long-time integrity of the simulation.

\subsection{Holland Thin-Wire Model}

We consider a closed wire antenna described by a smooth, closed curve $\Gamma$. The curve is discretized into $N_w$ piecewise-linear panels connecting consecutive vertices $\{\vX_q\}_{q=1}^{N_w}$, with $\vX_{N_w+1} = \vX_1$ for the closed wire. Panel~$q$ has length $\Delta X_q = |\vX_{q+1} - \vX_q|$, midpoint $\vX_{q+1/2} = (\vX_q + \vX_{q+1})/2$, and tangent direction $\vtau_q = (\vX_{q+1} - \vX_q)/\Delta X_q$.

The wire current $I_q^{n+1/2}$ is defined at panel midpoints and updated in sync with the magnetic field, while the voltage $V_q^n$ is defined at panel vertices and updated in sync with the electric field. This staggering mirrors the spatial staggering of the Yee grid and permits a fully explicit leapfrog scheme. The discrete telegrapher equations are
\begin{align}
L \frac{I_q^{n+1/2} - I_q^{n-1/2}}{\dt} &= -\frac{V_{q+1}^n - V_q^n}{\Delta X_q}
  + \Ical(\vE^n)_q - R I_q^n, \label{eq:Iupdate} \\
C \frac{V_q^{n+1} - V_q^n}{\dt} &= -\frac{I_q^{n+1/2} - I_{q-1}^{n+1/2}}{\Delta X_q},
  \label{eq:Vupdate}
\end{align}
in which $I_q^n = (I_q^{n+1/2} + I_q^{n-1/2})/2$ is the time-averaged current appearing in the resistive term. Throughout this paper we take the wire to be a perfect electrical conductor (PEC), so that $R = 0$. The operator $\Ical(\vE^n)_q$ interpolates the tangential electric field to panel~$q$. We will define the interpolation operator precisely, along with the current deposition operator, in the following section.

We remark that in the continuous setting the tangential electric field vanishes on a PEC surface. In Holland and Simpson's thin-wire formalism, however, $\Ical(\vE^n)_q$ represents the tangential field sampled at a relatively small distance from the wire surface that is large enough to be away from the singular near field of the wire, yet still well within a wavelength $\lambda$ of the antenna. This is the standard interpretation of the thin-wire model and is consistent with the slender-body approximation underlying the formalism~\cite{holland1981}.

The per-unit-length inductance $L$ and capacitance $C$ in~\eqref{eq:Iupdate}--\eqref{eq:Vupdate} are those of a thin wire of physical radius $a$ embedded in the regularized current density deposited onto the Yee grid. Following Holland and Simpson~\cite{holland1981} and B\'{e}renger~\cite{berenger2019},
\begin{equation}
L = \frac{\mu_0}{2\pi}\,\ln\!\left(\frac{d_{\mathrm{avg}}}{a}\right),
\qquad
C = \frac{\mu_0\,\epsilon_0}{L},
\label{eq:wire_LC}
\end{equation}
in which $d_{\mathrm{avg}}$ is the geometric-mean distance between the wire axis and the support of the regularized delta function $\delta_h$ on the grid, and the capacitance is fixed by the requirement that signals propagate along the wire at the free-space speed of light. The value of $d_{\mathrm{avg}}$ depends on the chosen $\delta_h$ and is computed once per choice of regularized delta function.

\subsection{Current Deposition and Field Interpolation}
\label{sec:operators}

The wire and grid are coupled through two operators. The \emph{current deposition} operator $\Scal$ maps the wire current to a source current density on the Yee grid, and the \emph{interpolation} operator $\Ical$ samples the tangential electric field from the grid back to the wire. Together they appear in the leapfrog updates~\eqref{eq:Hupdate}--\eqref{eq:Eupdate} and the telegrapher equations~\eqref{eq:Iupdate}--\eqref{eq:Vupdate}.

In the continuous setting, the current density associated with a wire carrying current $I(s)$ along the closed curve $\Gamma$ is
\begin{equation}
    \vJ(\vx) = \int_{\Gamma} I(s)\, \vtau(s)\,
    \delta(\vx - \vX(s))\, \mathrm{d}s,
    \label{eq:J_continuous}
\end{equation}
in which $\delta$ is the three-dimensional Dirac delta and $s$ is the arc-length parameter along the wire. This is a distributional vector field supported on $\Gamma$. The continuous interpolation operator evaluates the tangential electric field along the wire as
\begin{equation}
    \Ical(\vE)(s) = \int_{\R^3} \vE(\vx) \cdot \vtau(s)\,
    \delta(\vx - \vX(s))\, \mathrm{d}\vx,
    \label{eq:interp_continuous}
\end{equation}
which is simply the restriction of the tangential component of $\vE$ to $\Gamma$. In the continuous setting, $\Scal$ and $\Ical$ are formal adjoints of one another with respect to the $L^2(\R^3)$ inner product on fields and the $L^2(\Gamma)$ inner product on wire quantities.

The discrete operators are obtained by replacing $\delta$ with a regularized delta function $\delta_h$ and discretizing the wire as piecewise-linear panels with piecewise-constant currents. The deposited current density $\vJ_h = \Scal(I_q)$ at Yee edge $\vx_{i,j,k}$ is
\begin{equation}
    \vJ^{n+\frac{1}{2}}_h(\vx_{i,j,k}) = \sum_{q=1}^{N_w} I^{n+\frac{1}{2}}_q
    \int_0^{\Delta X_q} \vtau_{q}\,\delta_h(\vx_{i,j,k} - \vX_q(s))\, \mathrm{d}s.
    \label{eq:spread}
\end{equation}
In full generality, the interpolation operator is any linear functional of the form
\begin{equation}
    \Ical(\vE^n_h)_q = \sum_{i,j,k} w_{q,i,j,k}\, \vE^n_{i,j,k} \cdot \vtau_{q},
    \label{eq:interp_general}
\end{equation}
in which $w_{q,i,j,k}$ are interpolation weights to be chosen. The current deposition operator~\eqref{eq:spread} induces a natural choice of weights:
\begin{equation}
    w_{q,i,j,k} = \int_0^{\Delta X_q} \delta_h(\vx_{i,j,k} - \vX_q(s))\, \mathrm{d}s,
    \label{eq:interp_adjoint}
\end{equation}
and with this choice $\Ical$ is by construction the discrete adjoint of $\Scal$, satisfying the discrete power identity
\begin{equation}
\langle \vE,\, \Scal(I_q) \rangle_h = \langle I,\, \Ical(\vE) \rangle_w,
\label{eq:adjoint}
\end{equation}
in which $\langle \cdot, \cdot \rangle_h$ is the natural $\ell^2$ inner product on the Yee grid and
\begin{equation}
\langle a,\, b \rangle_w \;\equiv\; \sum_{q=1}^{N_w} a_q\, b_q\, \Delta X_q
\label{eq:wire_inner}
\end{equation}
is the corresponding $\ell^2$ inner product on the discretized wire. The left-hand side of~\eqref{eq:adjoint} is the power exchanged between the deposited current and the grid electric field, and the right-hand side is the power exchanged between the interpolated tangential field and the wire current. Without this identity the spatial discretization cannot be skew-symmetric, and energy can be spuriously created or destroyed by the discretization, undermining the conservation properties of the leapfrog scheme. For now, we leave both the choice of regularized delta function $\delta_h$ and the interpolation weights~\eqref{eq:interp_adjoint} unspecified. We specify their precise forms in Sections~\ref{sec:dist_curl} and~\ref{sec:adjoint} by the charge-conservation and adjointness requirements developed in the following section.

\section{Discrete Charge Conservation}
\label{sec:dist_curl}

\subsection{The Continuous Case}

We begin by examining the current density~\eqref{eq:J_continuous} in the continuous setting. Taking the divergence in the distributional sense and using the identity
\begin{equation}
\vtau(s) \cdot \nabla_\vx \delta(\vx - \vX(s)) = -\frac{\mathrm{d}}{\mathrm{d}s} \delta(\vx - \vX(s)),
\label{eq:tau_grad_delta}
\end{equation}
followed by integration by parts along the closed curve $\Gamma$, we obtain
\begin{equation}
\nabla \cdot \vJ
= \int_\Gamma \frac{\mathrm{d}I}{\mathrm{d}s}\, \delta(\vx - \vX(s))\,
\mathrm{d}s,
\label{eq:divJ_general}
\end{equation}
in which the boundary term vanishes since $\Gamma$ is closed. For a non-constant current, this is a distribution supported on $\Gamma$ with density $\mathrm{d}I/\mathrm{d}s$, consistent with the continuity equation $\frac{\partial\rho}{\partial t} + \nabla \cdot \vJ = 0$ via the telegrapher equation~\eqref{eq:Vupdate}. Moreover, \eqref{eq:divJ_general} vanishes if and only if $I$ is constant along $\Gamma$. In the distributional sense, if $\nabla\cdot\vJ\equiv 0$ then for any smooth compactly supported function $\phi$
\begin{equation}
0=\langle \nabla\cdot\vJ,\phi\rangle
=\int_\Gamma \frac{\mathrm{d}I}{\mathrm{d}s}(s)\,\phi(\vX(s))\,\mathrm{d}s,
\end{equation}
and since $\phi(\vX(\cdot))$ may range over all smooth functions on $\Gamma$, it follows that $\mathrm{d}I/\mathrm{d}s = 0$ and hence $I$ is constant. A constant current is therefore the unique current distribution on a closed wire that is divergence-free in the distributional sense, and correspondingly induces no charge accumulation anywhere in the domain. This observation motivates B\'{e}renger's criterion~\cite{berenger2019}, which is based on analyzing the deposition of a constant current density onto the grid.

\subsection{The Discrete Charge-Conservation Condition}

In the discrete setting, the current deposition operator $\Scal$ maps the wire current to a grid field $\vJ_h$. From~\eqref{eq:gauss_update}, Gauss's law is preserved if $\divh \vJ_h = 0$. The previous subsection shows that in the continuous setting, a constant current is the unique current distribution on a closed wire whose deposited field is divergence-free. The natural discrete analog is therefore to require that a constant wire current $I_q = I$ deposits a discretely divergence-free grid field:
\begin{equation}
\divh \vJ_h = 0 \quad \text{if } I_q = I \text{ for all } q.
\label{eq:berenger_cond}
\end{equation}
This is precisely B\'{e}renger's charge-conservation criterion~\cite{berenger2019}, which he derived through a geometric Q-cell argument without reference to the distributional nature of the source term. By the fundamental identity~\eqref{eq:div_curl_zero}, condition~\eqref{eq:berenger_cond} holds if the deposited current lies in the discrete curl subspace of the Helmholtz decomposition~\eqref{eq:helmholtz}. In the following section, we construct regularized delta functions $\delta_h$ that satisfy B\'{e}renger's criterion~\eqref{eq:berenger_cond} exactly.

\subsection{Kernel Compatibility with the Yee Staggering}

We now ask what structure the regularized delta function $\delta_h$ must have in order for condition~\eqref{eq:berenger_cond} to be satisfied. Writing $\delta_h$ as a tensor product of kernel functions
\begin{equation}
\delta_h(\vx) = \frac{1}{\dx\,\dy\,\dz}
\phi_x\!\Bigl(\frac{x}{\dx}\Bigr)
\phi_y\!\Bigl(\frac{y}{\dy}\Bigr)
\phi_z\!\Bigl(\frac{z}{\dz}\Bigr),
\label{eq:kernel_tensor}
\end{equation}
we examine the discrete divergence of a constant current ($I_q = I$ for each $q$) deposited on the Yee grid.

The $x$-component of the deposited current at Yee edge $(i+\frac{1}{2}, j, k)$ is
\begin{IEEEeqnarray}{rCl}
J_x^h\!\bigl(i\!+\!\tfrac{1}{2},j,k\bigr)
  &=& I \int \tau_x(s)\,
      \phi_x\!\Bigl(\tfrac{x_{i+1/2} - X_x(s)}{\dx}\Bigr) \eqieeebreak
  \eqieeecontinuequad
      \phi_y\!\Bigl(\tfrac{y_j - X_y(s)}{\dy}\Bigr)
      \phi_z\!\Bigl(\tfrac{z_{k} - X_z(s)}{\dz}\Bigr) \,\mathrm{d}s,
\IEEEyesnumber\label{eq:Jx_spread}
\end{IEEEeqnarray}
and analogous expressions hold for $J_y^h$ and $J_z^h$ with the appropriate staggering. The $x$-contribution to the discrete divergence at cell $(i,j,k)$ is the finite difference
\begin{equation}
(\nabla_h\cdot \vJ^h)_{x}\bigl(i,j,k\bigr)
:= \frac{1}{\dx}\Bigl(
J_x^h\!\bigl(i\!+\!\tfrac{1}{2},j,k\bigr)
- J_x^h\!\bigl(i\!-\!\tfrac{1}{2},j,k\bigr)\Bigr),
\label{eq:Jx_div_def}
\end{equation}
which introduces the centered difference
\begin{equation}
\frac{1}{\dx}\left[
\phi_x\!\Bigl(\tfrac{x_{i+1/2} - X_x(s)}{\dx}\Bigr)
- \phi_x\!\Bigl(\tfrac{x_{i-1/2} - X_x(s)}{\dx}\Bigr)
\right].
\label{eq:Jx_diff}
\end{equation}
If the function $\phi_x$ is such that its staggered difference is equal to the true continuous derivative of another function $\psi_x$, then we can apply the chain rule to the difference in~\eqref{eq:Jx_diff} to get
\begin{equation}
\frac{1}{\dx}\left[
\phi_x\!\Bigl(r + \tfrac{1}{2}\Bigr) - \phi_x\!\Bigl(r - \tfrac{1}{2}\Bigr)
\right]
= \frac{\mathrm{d}}{\mathrm{d}r}\psi_x(r).
\label{eq:kernel_cond}
\end{equation}
Setting $r=(x_i-X_x(s))/\dx$, this converts the $x$-contribution to the discrete divergence into a term involving $\frac{\partial}{\partial X_x}\psi_x$. Applying the same construction to the $y$- and $z$-differences, the total discrete divergence takes the form
\begin{IEEEeqnarray}{rCl}
\nabla_h\cdot \vJ^h_{i,j,k}
  &=& I\int
      \frac{\mathrm{d}}{\mathrm{d}s}\Bigl[
      \psi_x\!\Bigl(\tfrac{x_i-X_x(s)}{\dx}\Bigr) \eqieeebreak
  \eqieeecontinuequad
      \psi_y\!\Bigl(\tfrac{y_j-X_y(s)}{\dy}\Bigr)
      \psi_z\!\Bigl(\tfrac{z_k-X_z(s)}{\dz}\Bigr)
      \Bigr]\,\mathrm{d}s.
\IEEEyesnumber\label{eq:perfect_diff}
\end{IEEEeqnarray}
Since $\Gamma$ is closed, the integral of a total $s$-derivative vanishes, and therefore $\nabla_h\cdot \vJ^h(i,j,k) = 0$, confirming~\eqref{eq:berenger_cond}.

Condition~\eqref{eq:kernel_cond} has a natural realization in terms of cardinal B-splines
$\mathrm{BS}_n$, defined by the recursive convolution~\cite{deboor1991}
\begin{equation}
\mathrm{BS}_{n+1}(r) = (\mathrm{BS}_n * \mathrm{BS}_0)(r),
\label{eq:bspline_recursion}
\end{equation}
in which $\mathrm{BS}_0$ is the indicator function of $[-\tfrac{1}{2}, \tfrac{1}{2})$. Differentiating~\eqref{eq:bspline_recursion} and applying the fundamental theorem of calculus yields
\begin{equation}
\frac{\mathrm{d}}{\mathrm{d}r}\mathrm{BS}_{n+1}(r)
=
\mathrm{BS}_n\!\bigl(r + \tfrac{1}{2}\bigr)
- \mathrm{BS}_n\!\bigl(r - \tfrac{1}{2}\bigr),
\label{eq:bspline_id}
\end{equation}
so taking $\phi_\alpha = \mathrm{BS}_n$ satisfies~\eqref{eq:kernel_cond} with $\psi_\alpha = \mathrm{BS}_{n+1}$. The node-aligned kernels for each component $J_\alpha$ are then also taken to be $\mathrm{BS}_{n+1}$, so that all three factors in~\eqref{eq:perfect_diff} are of the same type and the chain rule identity applies uniformly. The resulting family of regularized delta functions is denoted $\delta_h^{(n)}$, in which the superscript indexes the order of the B-spline kernel in the staggered direction, and is defined by
\begin{equation}
\delta_h^{(n)}(\vx - \vX) = \frac{1}{\dx\,\dy\,\dz}\, \mathrm{BS}_n\!\left(\frac{x_\alpha - X_\alpha}{\Delta x_\alpha}\right) \prod_{\beta \neq \alpha} \mathrm{BS}_{n+1}\!\left(\frac{x_\beta - X_\beta}{\Delta x_\beta}\right),
\label{eq:delta_h_family}
\end{equation}
in which $\alpha \in \{x,y,z\}$ is the Cartesian direction along which the deposited Yee component is staggered ($\alpha = x$ for $J_x$, etc.) and $\Delta x_\alpha$ denotes the grid spacing in direction $\alpha$. In Section~\ref{sec:results} we test $\delta_h^{(0)}$, $\delta_h^{(2)}$, and $\delta_h^{(4)}$, corresponding to staggered orders $\mathrm{BS}_0$, $\mathrm{BS}_2$, $\mathrm{BS}_4$ paired respectively with $\mathrm{BS}_1$, $\mathrm{BS}_3$, $\mathrm{BS}_5$ in the two node-aligned directions. We remark that B-splines are not the only admissible kernels: any $\phi_\alpha$ satisfying~\eqref{eq:kernel_cond} for some $\psi_\alpha$ would suffice. By contrast, the isotropic choice $\delta_h^{\mathrm{iso}}$, with the piecewise linear $\mathrm{BS}_1$ in all three directions, does not satisfy~\eqref{eq:kernel_cond} and produces discrete divergence errors regardless of the accuracy of the quadrature approximation used.

Similar composite B-spline constructions have appeared previously in compatible discretizations of Maxwell's equations in the context of isogeometric analysis~\cite{buffa2010,simpson2018}, in which the same composite pairing arises from the requirement that the discrete de Rham complex be exact. A related construction was used in~\cite{roy-chowdhury2024} to interpolate discretely curl-free vector fields on a staggered grid to continuously curl-free data, exploiting the same algebraic compatibility between B-spline orders and grid staggering developed here. We also note that the lowest-order member $\delta_h^{(0)}$ coincides with the lowest-order Whitney 1-form on the edges of a cube~\cite{bossavit2000,glasser2023}, which provides a finite element analog of the Yee grid discretization in the language of differential forms.

We note that the lowest-order member $\delta_h^{(0)}$ has the same form as the regularized delta function underlying the deposition rule of Guiffaut, Reineix, and Pecqueux~\cite{guiffaut2010,guiffaut2012}, with their endpoint-averaging step replaced here by exact Gauss--Legendre quadrature of the line integral~\eqref{eq:perfect_diff}. B\'{e}renger's Q-cell construction~\cite{berenger2019}, although motivated geometrically rather than algebraically, also fits naturally within this lowest-order framework: the Q-cell facet pierced by the wire selects a piecewise-constant assignment in the staggered, normal direction, and the linear interpolation of the wire current to the surrounding $J$-nodes plays the role of the bilinear transverse weighting in~\cite{guiffaut2012}.

B\'{e}renger does not, however, define an interpolation operator adjoint to his current deposition. He observes that the resulting asymmetry leaves the leapfrog updates stable in his simulations, and on that basis does not pursue an adjoint construction. However, we note that the adjoint property is less a stability requirement than a statement of energy conservation: the discrete power identity it implies (Section~\ref{sec:adjoint}) is what ensures the spatial discretization remains skew-symmetric so that the leapfrog update does not dissipate or create any energy.

\section{Exact Panel Quadrature}
\label{sec:kernels}

Maintaining charge conservation requires both that the kernel pair $(\phi_\alpha, \psi_\alpha)$ satisfy the compatibility condition~\eqref{eq:kernel_cond} and that the line integral in~\eqref{eq:perfect_diff} be evaluated exactly. The kernel pair guarantees that the integrand has the form of a total derivative along the wire, so its value vanishes on a closed loop if the integral is evaluated exactly. Without exact evaluation, the total-derivative identity breaks at the discrete level and a residual charge proportional to the quadrature error is deposited.

The cardinal B-spline family $\delta_h^{(n)}$, with $\mathrm{BS}_n$ in the staggered direction and $\mathrm{BS}_{n+1}$ in the two node-aligned directions, is one such compatible pair; increasing $n$ produces smoother current distributions on the grid at the cost of wider stencils. Exact panel quadrature is possible here precisely because the B-spline kernels are piecewise polynomial. We split each panel at the points where the wire crosses a transition between the polynomial pieces of one of the kernels. On each resulting subpanel, the kernel factors have fixed polynomial representations, so the integrand is polynomial in the arc-length parameter $s$. Composite Gauss--Legendre quadrature with subinterval breakpoints at these crossings, and with enough Gauss points per subinterval to integrate the polynomial pieces exactly, therefore evaluates the panel integral to machine precision. Other compatible pairs are admissible, but for kernels that are not piecewise polynomial, an exact quadrature scheme is generally not available. The same caveat applies to the wire representation: the integrand's polynomial form in $s$ relies on the piecewise-linear panel parameterization $\vX_q(s) = \vX_q + s\,\vtau_q$. A higher-order wire representation, such as a cubic-spline parameterization, would make $\vX_q(s)$ a non-polynomial function of arc length, and the integrand would acquire a non-polynomial dependence on $s$ that an exact quadrature would have to address.

Due to the staggering of the Yee grid, for the $\delta_h^{(n)}$ family, even $n$ places the breakpoints at integer grid-plane crossings of the wire, while odd $n$ places them at half-integer crossings. The quadrature procedure is otherwise the same in either case.

\section{Zero Loop EMF for Discrete Gradient Fields}
\label{sec:adjoint}

In Section~\ref{sec:operators} we defined the interpolation operator $\Ical$ as the discrete adjoint of the current deposition operator $\Scal$ by choosing the interpolation weights~\eqref{eq:interp_adjoint} to equal the current deposition weights, and stated the resulting power identity~\eqref{eq:adjoint}. Combined with the leapfrog updates of Section~\ref{sec:background}, the power identity~\eqref{eq:adjoint} ensures that the spatial discretization remains skew-symmetric and that no spurious energy is created or destroyed in the coupled FDTD--thin-wire system. A complete proof of this energy-conservation property under periodic boundary conditions, together with the resulting CFL-type time-step restriction, is given in the accompanying supplementary material.

The adjoint construction gives $\Ical$ the following discrete compatibility property: a discrete gradient field produces zero electromotive force around a closed wire loop. That is,
\begin{equation}
\vE_h = \nabla_h \phi
\quad\implies\quad
\langle \mathbf{1},\, \Ical(\vE_h) \rangle_w = 0,
\label{eq:emf_zero}
\end{equation}
in which $\mathbf{1}$ denotes the constant unit current on the wire degrees of freedom. Thus
\[
\langle \mathbf{1},\, \Ical(\vE_h) \rangle_w
=
\sum_q \Ical(\vE_h)_q \,\Delta X_q
\]
is the discrete electromotive force around the loop.

To prove \eqref{eq:emf_zero}, let $I_q = I$ be any constant current on the closed wire. By the charge-conservation property of the spreading operator,
$\Scal(I)$ is discretely divergence-free. Hence, by the discrete Helmholtz decomposition \eqref{eq:helmholtz}, it is equal to the discrete curl of some vector potential. Since discrete gradient fields are $\ell^2$ orthogonal to discrete curls on the Yee grid, we have
\[
\langle \nabla_h \phi,\, \Scal(I)\rangle_h = 0.
\]
Therefore, via the adjoint identity \eqref{eq:adjoint}, we obtain
\begin{equation}
I\,\langle \mathbf{1},\, \Ical(\nabla_h \phi) \rangle_w
=
\langle \nabla_h \phi,\, \Scal(I)\rangle_h
=
0.
\label{eq:emf_zero_proof}
\end{equation}
Since this holds for any constant $I$, the discrete loop EMF must vanish, establishing \eqref{eq:emf_zero}. This is the discrete analog of the fact that an electrostatic field has zero circulation around a closed loop, which is the magnetostatic version of Faraday's law.

The interpolation operator therefore preserves the property that discretely irrotational electric fields produce no net electromotive force around a closed loop. This compatibility can fail for isotropic tensor-product regularized delta functions, because the associated deposition operator generally does not satisfy~\eqref{eq:berenger_cond}. In that case, a constant closed-loop current need not deposit to a discretely divergence-free grid current, and the orthogonality argument above breaks down. Our construction achieves charge-conserving deposition and zero-EMF interpolation for gradient fields simultaneously, with the adjoint identity providing the link between them.

\section{Numerical Results}
\label{sec:results}

We test these operators in three computational experiments: the input impedance and current of a center-fed dipole, a circular loop, and a square loop. All experiments are run on a uniform three-dimensional Yee grid with $N = 32$ cells per unit length ($\dx = \dy = \dz = h = 1/N$), bounded by a 32-cell CFS-PML absorbing layer~\cite{gedney2011}. The wire is discretized into uniform panels of length $\ds \approx h$. Each panel has radius $a = h / 10$, and the per-unit-length parameters $L$ and $C$ are evaluated from~\eqref{eq:wire_LC} with $d_{\mathrm{avg}}$ taken from the chosen $\delta_h^{(n)}$. The time step satisfies the CFL condition $\dt = h / (2\sqrt{3}\,c_0)$. We note that this is the CFL condition associated with the FDTD method alone~\cite{taflove2005}. This is sufficient for stability because we assume that the wire radius $a$ is relatively small and consistent with the thin-wire formalism. If the wire radius is larger, the time-step size will need to be reduced to maintain stability~\cite{berenger2021,ledfelt2002}.

In every experiment we compare four regularized delta functions: the composite charge-conserving members $\delta_h^{(0)}$, $\delta_h^{(2)}$, and $\delta_h^{(4)}$ of the family defined in~\eqref{eq:delta_h_family}, and the piecewise-linear isotropic $\delta_h^{\mathrm{iso}}$ as a control that does not satisfy condition~\eqref{eq:kernel_cond}. We evaluate three orientations of each wire with respect to the grid: an axis-aligned configuration (loop normal $\hat{n}=\hat{z}$), an in-plane face-diagonal configuration ($\hat{n}=(\hat{x}+\hat{y})/\sqrt{2}$), and a body-diagonal configuration ($\hat{n}=(\hat{x}+\hat{y}+\hat{z})/\sqrt{3}$). These three scenarios run from grid-aligned to fully off-axis. Since the antenna response is rotationally invariant in free space, the curves should ideally be identical.

\subsection{Linear Dipole Antenna}
\label{sec:dipole_imp}

We consider a center-fed straight dipole of length $L_{\mathrm{dipole}} = 0.5\,\mathrm{m}$ and three orientations relative to the grid. The wire is open at both ends, so the current vanishes at the tips. A Gaussian voltage pulse is applied across a delta-gap at the center panel, with pulse width chosen so that the spectrum extends to $L_{\mathrm{dipole}}/\lambda \approx 2$, covering the first four resonances and antiresonances of the dipole. From the recorded gap voltage and current, we compute the input impedance and report the input resistance $R_{\mathrm{in}}(f)$ and reactance $X_{\mathrm{in}}(f)$ as functions of the electrical length $L_{\mathrm{dipole}}/\lambda$. In the thin-wire limit, the first resonance occurs at $L_{\mathrm{dipole}}/\lambda \approx 1/2$, with input resistance $73\,\Omega$ and reactance $42.5\,\Omega$~\cite{balanis2005}.

\begin{figure}[t]
\centering
\subfloat[$\delta_h^{(0)}$]{%
  \includegraphics[width=\resultfigwidth]{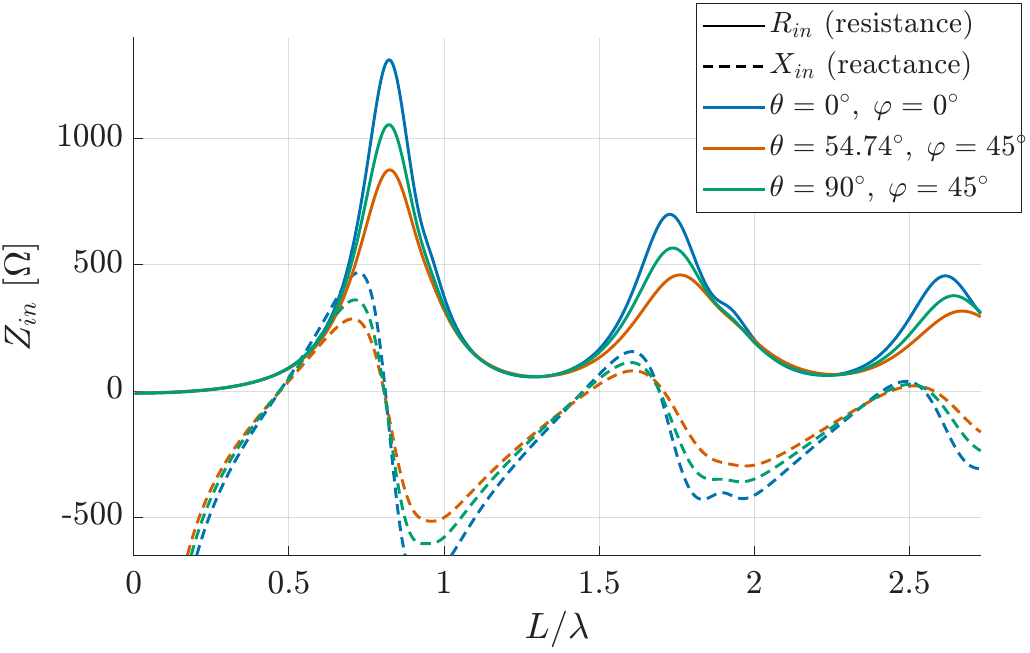}%
  \label{fig:dipole_Zin_BS0BS1}}\resultfigcolsep
\subfloat[$\delta_h^{(2)}$]{%
  \includegraphics[width=\resultfigwidth]{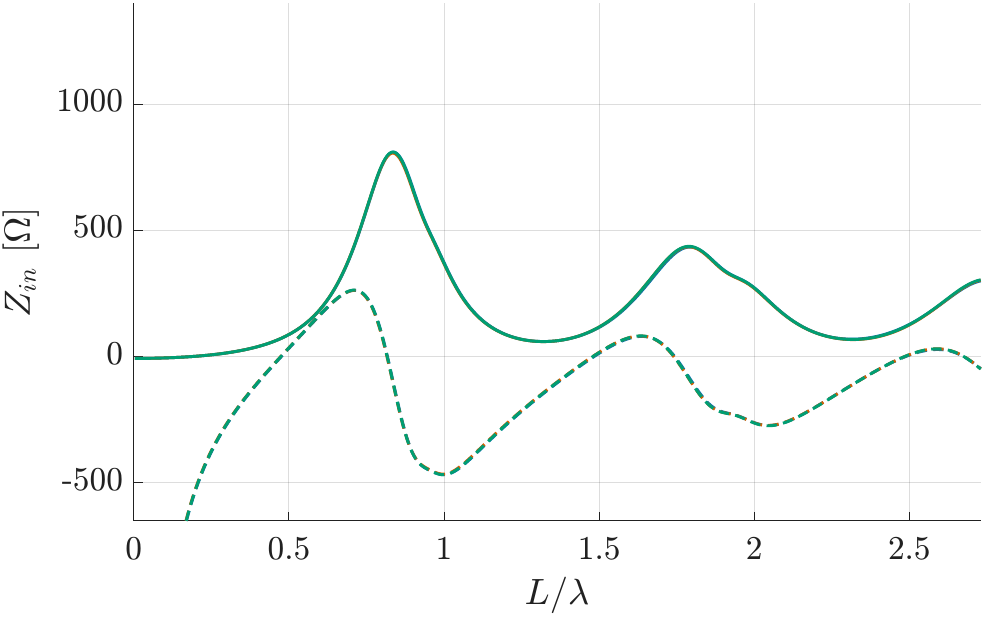}%
  \label{fig:dipole_Zin_BS2BS3}}\resultfigrowsep
\subfloat[$\delta_h^{(4)}$]{%
  \includegraphics[width=\resultfigwidth]{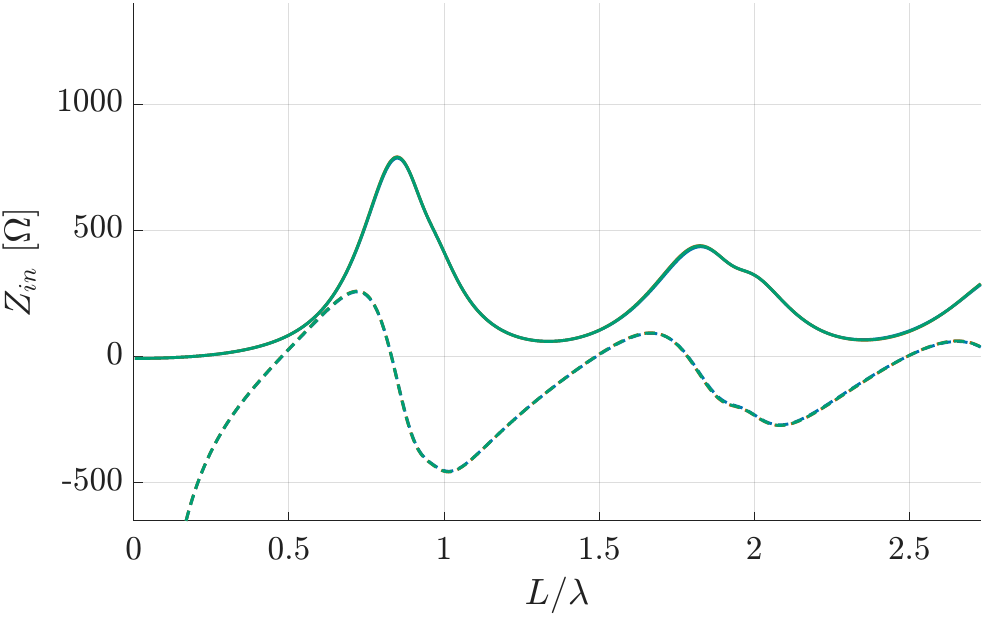}%
  \label{fig:dipole_Zin_BS4BS5}}\resultfigcolsep
\subfloat[$\delta_h^{\mathrm{iso}}$]{%
  \includegraphics[width=\resultfigwidth]{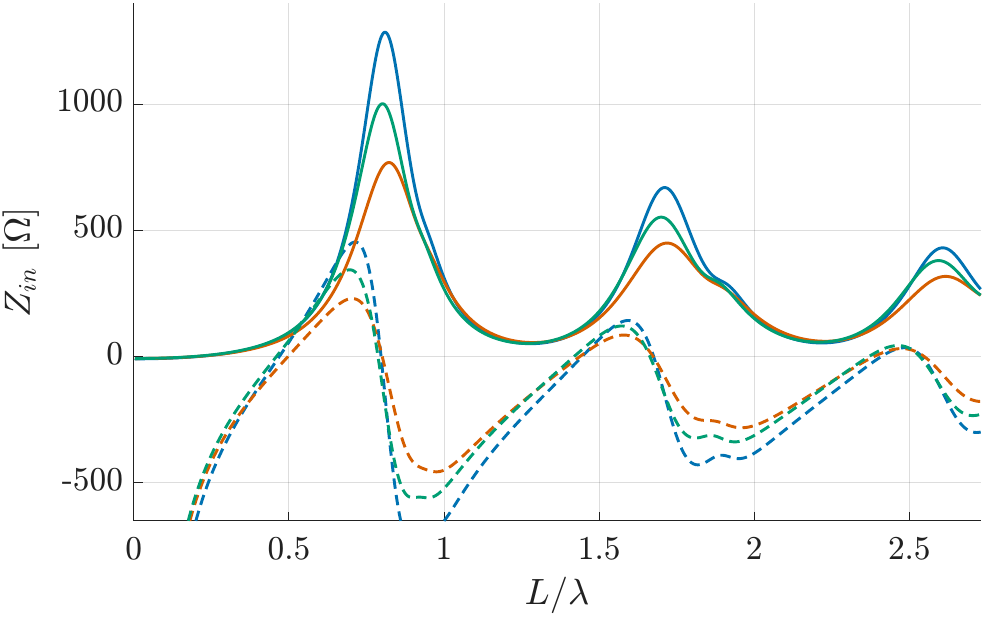}%
  \label{fig:dipole_Zin_iso}}
\caption{Input impedance of the center-fed dipole versus dipole length in wavelengths $L_{\mathrm{dipole}}/\lambda$, with resistance $R_{\mathrm{in}}$ as solid lines and reactance $X_{\mathrm{in}}$ as dashed lines. Subpanels (a)--(d) use the four regularized delta functions $\delta_h^{(0)}$, $\delta_h^{(2)}$, $\delta_h^{(4)}$, and the isotropic $\delta_h^{\mathrm{iso}}$. Within each subpanel, three orientation curves are overlaid: axis-aligned $\hat{\xi}=\hat{z}$, face-diagonal $\hat{\xi}=(\hat{x}+\hat{y})/\sqrt{2}$, and body-diagonal $\hat{\xi}=(\hat{x}+\hat{y}+\hat{z})/\sqrt{3}$.}
\label{fig:dipole_Zin}
\end{figure}

Fig.~\ref{fig:dipole_Zin} reports $R_{\mathrm{in}}$ and $X_{\mathrm{in}}$ for the three orientations and four regularized delta functions. The dipole is open at both tips, so the homogeneous boundary conditions on the wire current preclude a constant null mode, and the persistent low-frequency oscillation observed on closed loops does not arise here. For this reason, we report only the input impedance. All four regularized delta functions resolve the resonance pattern of a thin finite-length dipole, with the half-wave resonance near $L_{\mathrm{dipole}}/\lambda = 1/2$ landing on the canonical values $R_{\mathrm{in}} \approx 73\,\Omega$ and $X_{\mathrm{in}} \approx 42.5\,\Omega$~\cite{balanis2005}. The kernels differ in how cleanly the three orientation curves collapse on top of one another. The composite $\delta_h^{(2)}$ and $\delta_h^{(4)}$ regularized delta functions produce orientation-independent curves to within plotting tolerance. The $\delta_h^{(0)}$ and $\delta_h^{\mathrm{iso}}$ regularized delta functions show some spread in the impedance values, especially near the antiresonant modes of the dipole antenna. The sensitivity to orientation of $\delta_h^{(0)}$ is consistent with the discontinuous nature of the $\mathrm{BS}_0$ kernel. In addition, $\mathrm{BS}_0$ does not satisfy the discrete first-moment condition~\cite{engquist2005,hosseini2016}, so on sufficiently smooth fields $\mathrm{BS}_0$ gives only first-order accurate interpolation, while the other three kernels tested are formally second-order accurate. For a zero-radius wire model, however, this accuracy statement should be interpreted as a formal consistency result. Line-supported charge and current sources produce fields that are singular on the wire, and that singularity limits the achievable interpolation accuracy regardless of the kernel's formal order. An analogous subtlety arises in fluid--structure interaction, in which the Stokeslet (the Green's function for Stokes flow) admits only first-order accurate interpolation even if the regularized delta function achieves higher order on smooth fields~\cite{mori2008,liu2012,griffith2005}.

\resultfloatbarrier
\subsection{Circular Loop Antenna}
\label{sec:circle_imp}

The second experiment uses a circular loop antenna of radius $0.5\,\mathrm{m}$, driven by a differentiated Gaussian voltage pulse across a delta-gap at a single panel. The differentiated Gaussian carries no DC content, which is important for closed-loop antennas because the periodic boundary conditions on the current admit a constant mode (Section~\ref{sec:dist_curl}). On a PEC wire, a DC driving voltage would inject a circulating current that persists indefinitely after the pulse ends. With the differentiated source, any low-frequency residual current in the simulation is therefore attributable to the discretization rather than the excitation. The pulse width is chosen so that the spectrum extends to $C/\lambda \approx 2$, covering the first three resonances of the loop, with $C$ denoting the loop circumference. As with the dipole, we report the input resistance and reactance as functions of the
dimensionless circumference $C/\lambda$. We also monitor the time-domain ringdown of the gap current.

\begin{figure}[t]
\centering
\subfloat[$\delta_h^{(0)}$]{%
  \includegraphics[width=\resultfigwidth]{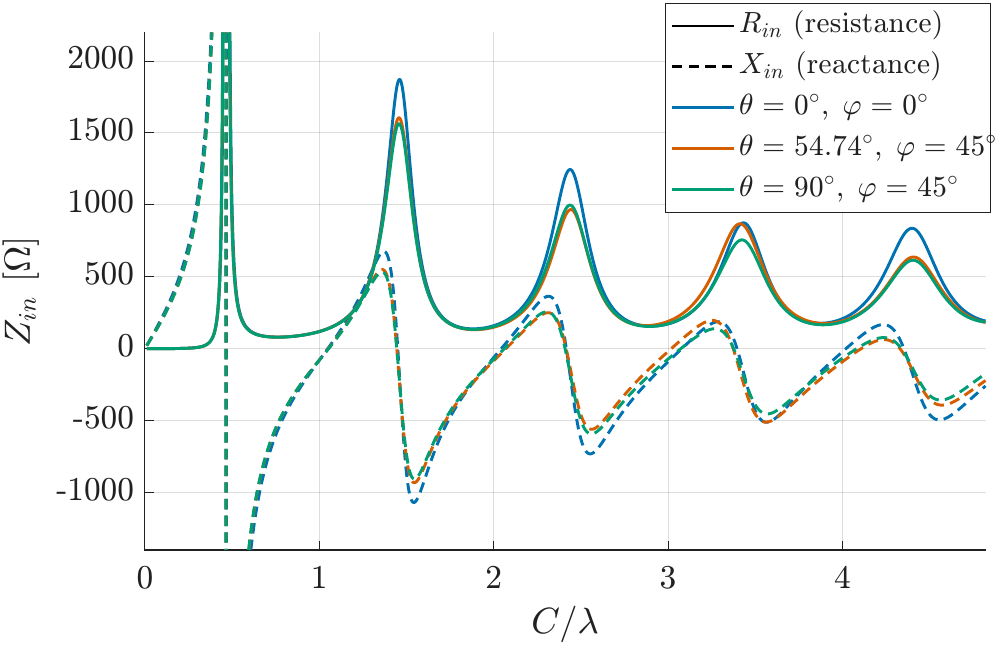}%
  \label{fig:circle_Zin_BS0BS1}}\resultfigcolsep
\subfloat[$\delta_h^{(2)}$]{%
  \includegraphics[width=\resultfigwidth]{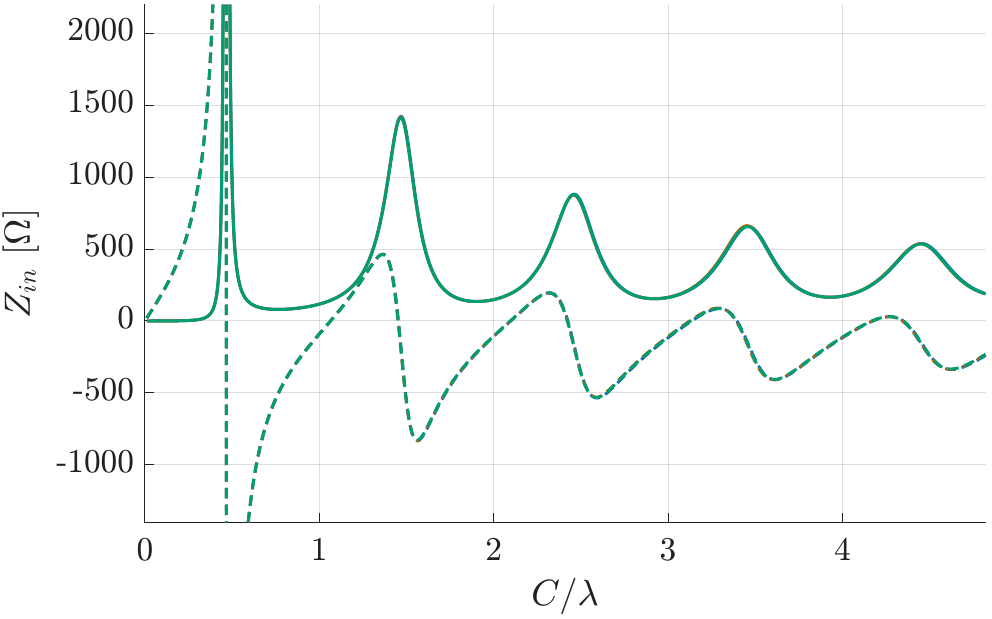}%
  \label{fig:circle_Zin_BS2BS3}}\resultfigrowsep
\subfloat[$\delta_h^{(4)}$]{%
  \includegraphics[width=\resultfigwidth]{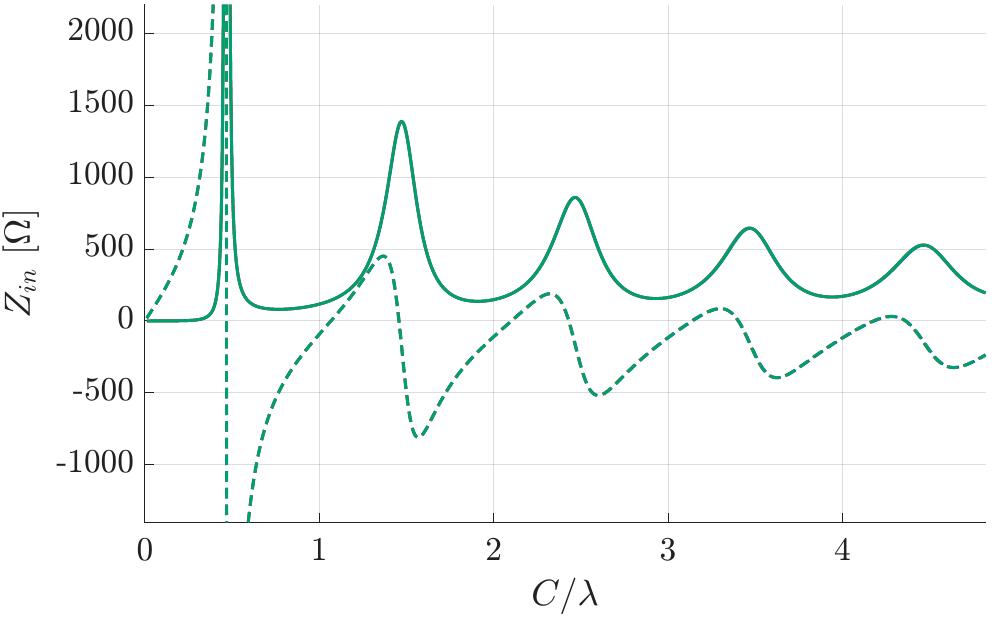}%
  \label{fig:circle_Zin_BS4BS5}}\resultfigcolsep
\subfloat[$\delta_h^{\mathrm{iso}}$]{%
  \includegraphics[width=\resultfigwidth]{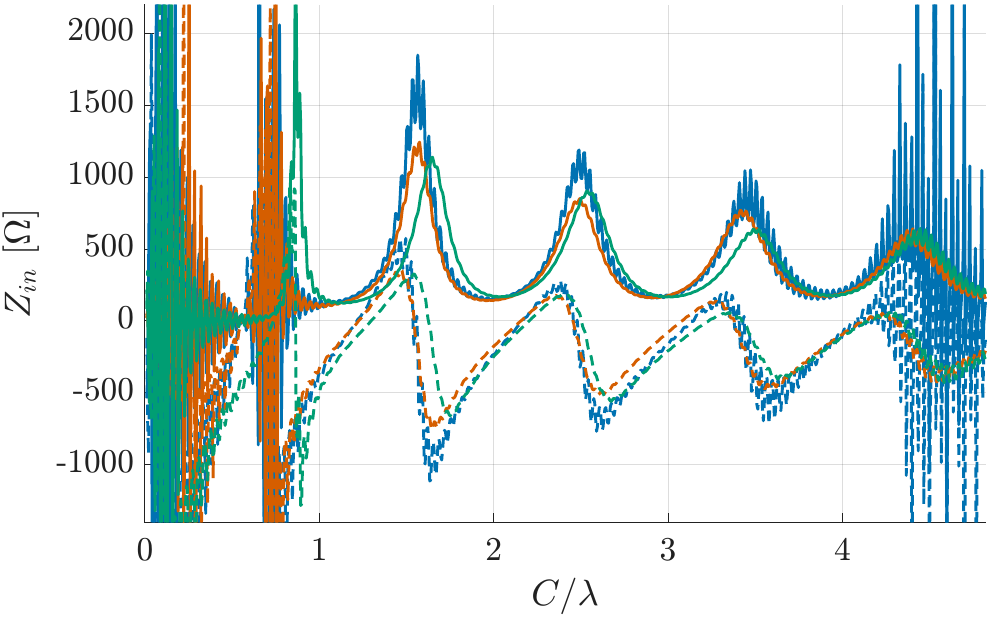}%
  \label{fig:circle_Zin_iso}}
\caption{Input impedance of the circular loop versus circumference in wavelengths $C/\lambda$, with resistance $R_{\mathrm{in}}$ as solid lines and reactance $X_{\mathrm{in}}$ as dashed lines. Subpanels (a)--(d) use the four regularized delta functions $\delta_h^{(0)}$, $\delta_h^{(2)}$, $\delta_h^{(4)}$, and the isotropic $\delta_h^{\mathrm{iso}}$. Within each subpanel, three orientation curves are overlaid: axis-aligned $\hat{n}=\hat{z}$, face-diagonal $\hat{n}=(\hat{x}+\hat{y})/\sqrt{2}$, and body-diagonal $\hat{n}=(\hat{x}+\hat{y}+\hat{z})/\sqrt{3}$.}
\label{fig:circle_Zin}
\end{figure}

\begin{figure}[t]
\centering
\subfloat[$\delta_h^{(0)}$]{%
  \includegraphics[width=\resultfigwidth]{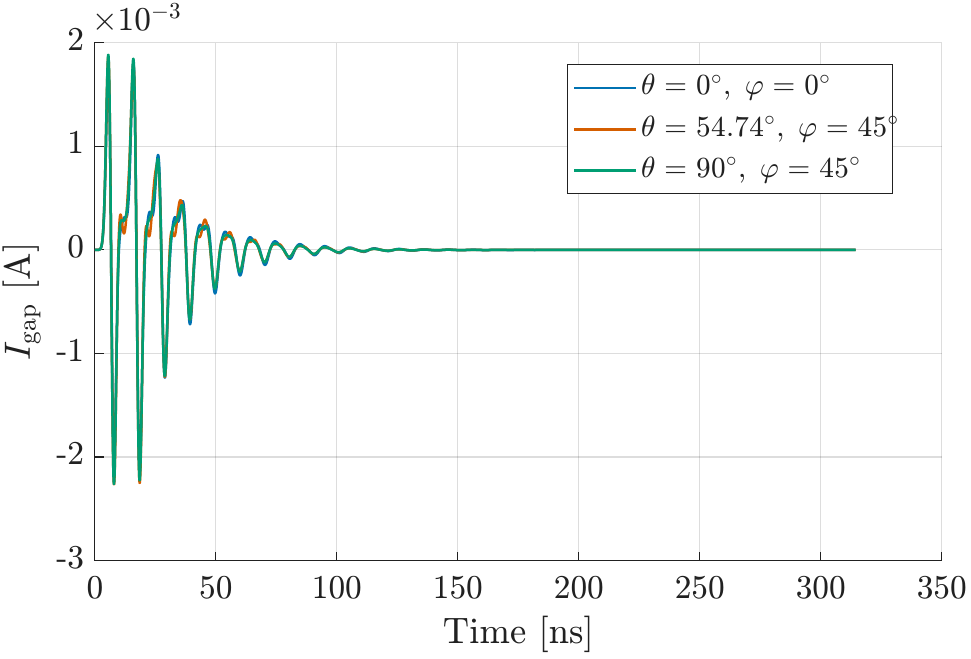}%
  \label{fig:circle_Igap_t_BS0BS1}}\resultfigcolsep
\subfloat[$\delta_h^{(2)}$]{%
  \includegraphics[width=\resultfigwidth]{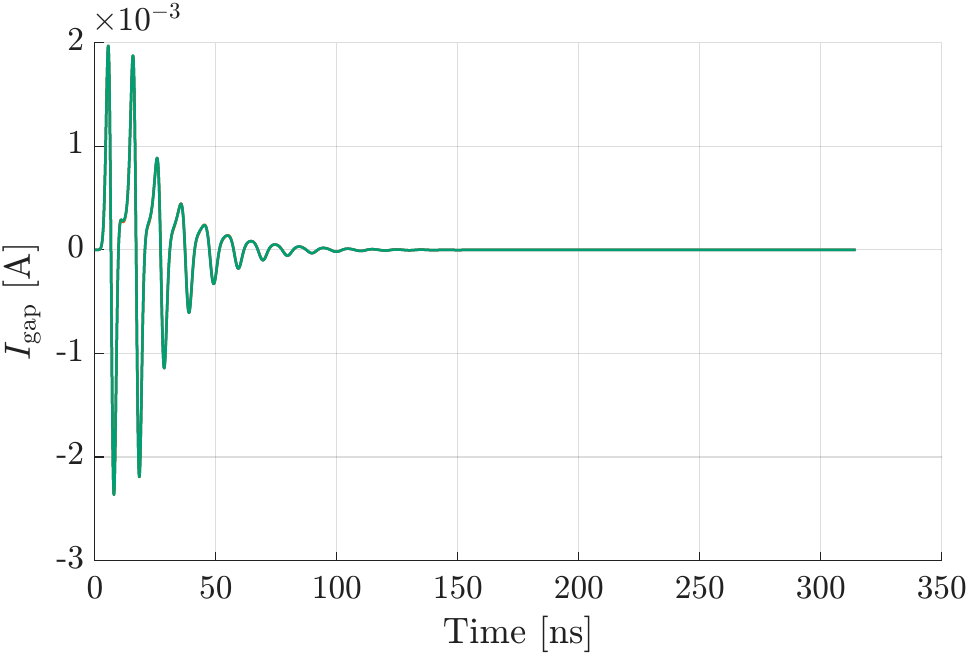}%
  \label{fig:circle_Igap_t_BS2BS3}}\resultfigrowsep
\subfloat[$\delta_h^{(4)}$]{%
  \includegraphics[width=\resultfigwidth]{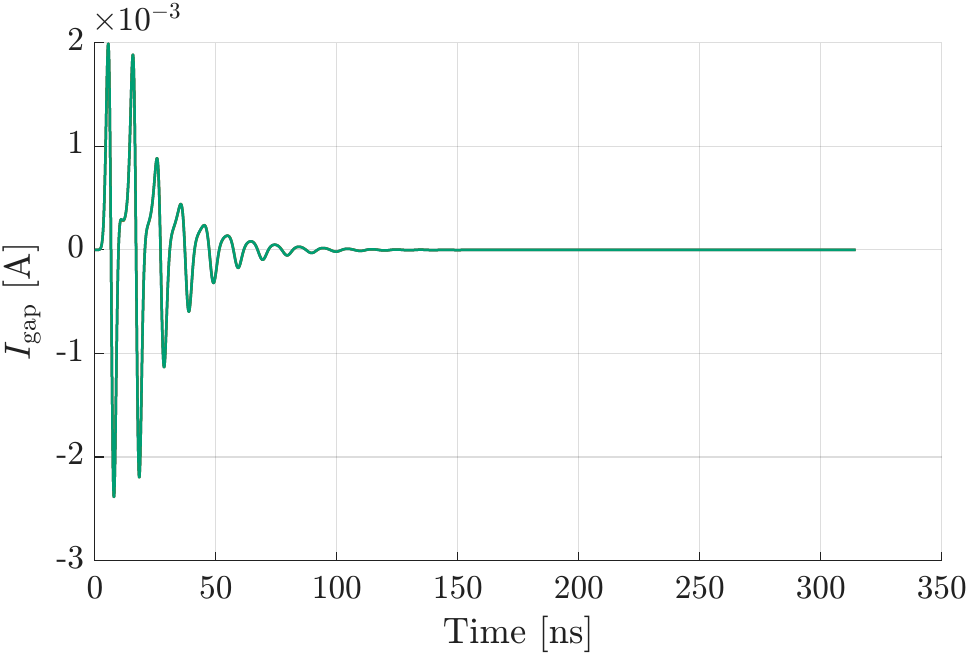}%
  \label{fig:circle_Igap_t_BS4BS5}}\resultfigcolsep
\subfloat[$\delta_h^{\mathrm{iso}}$]{%
  \includegraphics[width=\resultfigwidth]{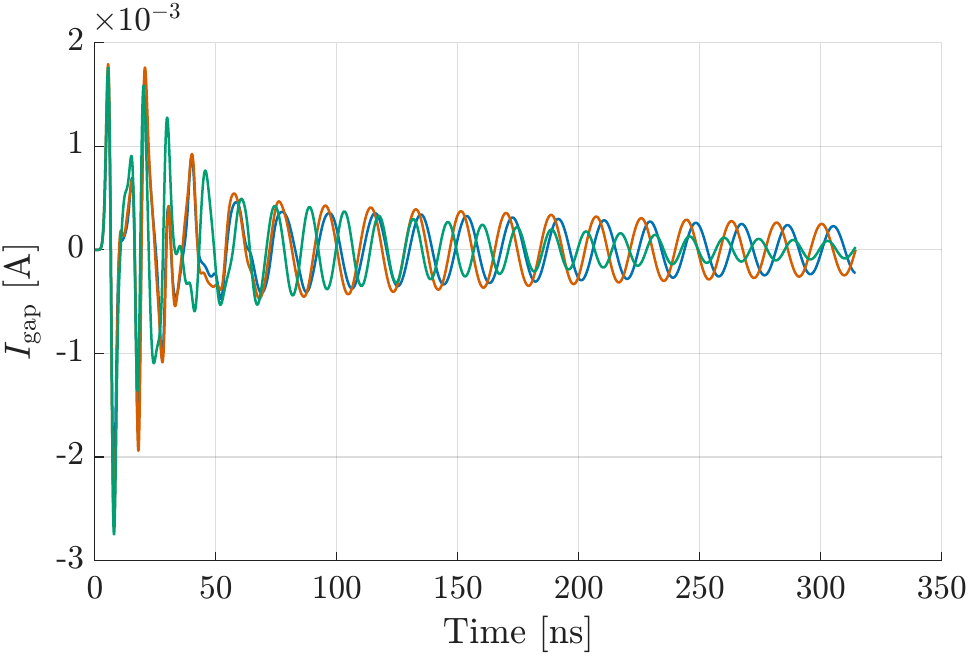}%
  \label{fig:circle_Igap_t_iso}}
\caption{Gap-current time history $I_{\mathrm{gap}}(t)$ at the feed of the circular loop. Subpanels (a)--(d) use the four regularized delta functions $\delta_h^{(0)}$, $\delta_h^{(2)}$, $\delta_h^{(4)}$, and the isotropic $\delta_h^{\mathrm{iso}}$, with the three orientation curves overlaid in each.}
\label{fig:circle_Igap_t}
\end{figure}

The circular loop antenna is well characterized analytically since the telegrapher equations on the loop can be solved more readily using Fourier series~\cite{storer1956,wu1962a,balanis2005,mckinley2013}, yielding closed-form expressions for the input impedance. The resonances of the loop antenna occur near $\lambda \approx C/m$ and the antiresonances near
$\lambda \approx 2C/m$, with $m$ a positive integer. Fig.~\ref{fig:circle_Zin} shows the input impedance and Fig.~\ref{fig:circle_Igap_t} the gap-current time histories, for three loop orientations and four kernel choices.

All composite regularized delta functions accurately reproduce the input impedance of the circular loop antenna, with resonances and antiresonances at their theoretically predicted locations. The gap current $I_{\mathrm{gap}}(t)$ decays to near zero within approximately $125\,\mathrm{ns}$, and the results are consistent across all three orientations, with minor variations for the least regular composite kernel $\delta_h^{(0)}$.

The isotropic regularized delta function $\delta_h^{\mathrm{iso}}$, by contrast, fails to reproduce the correct impedance characteristics. The impedance plot is dominated by noise, and the gap current develops a persistent low-frequency oscillation~\cite{edelvik2003,berenger2019,guiffaut2010,ledfelt2002} that shows no sign of decay over the duration of the simulation. This undamped oscillation renders the impedance nonphysical.

\resultfloatbarrier
\subsection{Square Loop Antenna}
\label{sec:square_imp}

In the final experiment, we replace the circular loop with a square loop antenna of side length $1\,\mathrm{m}$, fed at the midpoint of one side by the same differentiated Gaussian pulse used in Section~\ref{sec:circle_imp}. The square loop is well characterized and behaves similarly to the circular loop~\cite{balanis2005,stutzman2013,king1956}: resonances occur near $\lambda \approx C/m$ and antiresonances near $\lambda \approx 2C/m$, with $m$ a positive integer and $C$ the perimeter. For small $C/\lambda$ the reactance is predominantly inductive, becoming capacitive at larger values before alternating between the two. This contrasts with the linear dipole, whose reactance begins capacitive before becoming inductive and then alternating between the two, as shown in Fig.~\ref{fig:dipole_Zin}.

\begin{figure}[t]
\centering
\subfloat[$\delta_h^{(0)}$]{%
  \includegraphics[width=\resultfigwidth]{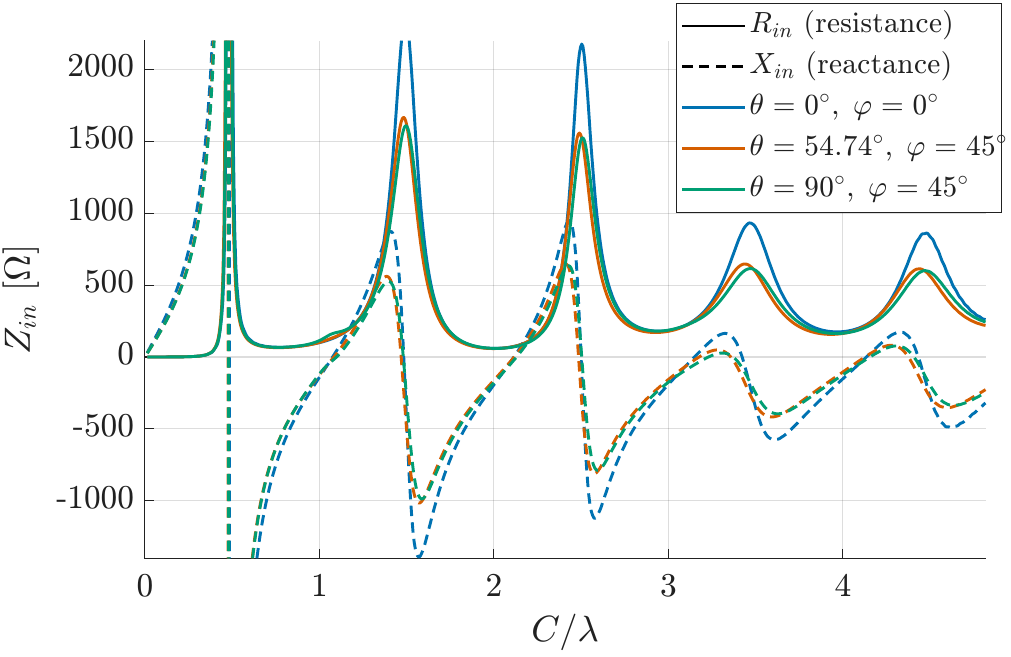}%
  \label{fig:square_Zin_BS0BS1}}\resultfigcolsep
\subfloat[$\delta_h^{(2)}$]{%
  \includegraphics[width=\resultfigwidth]{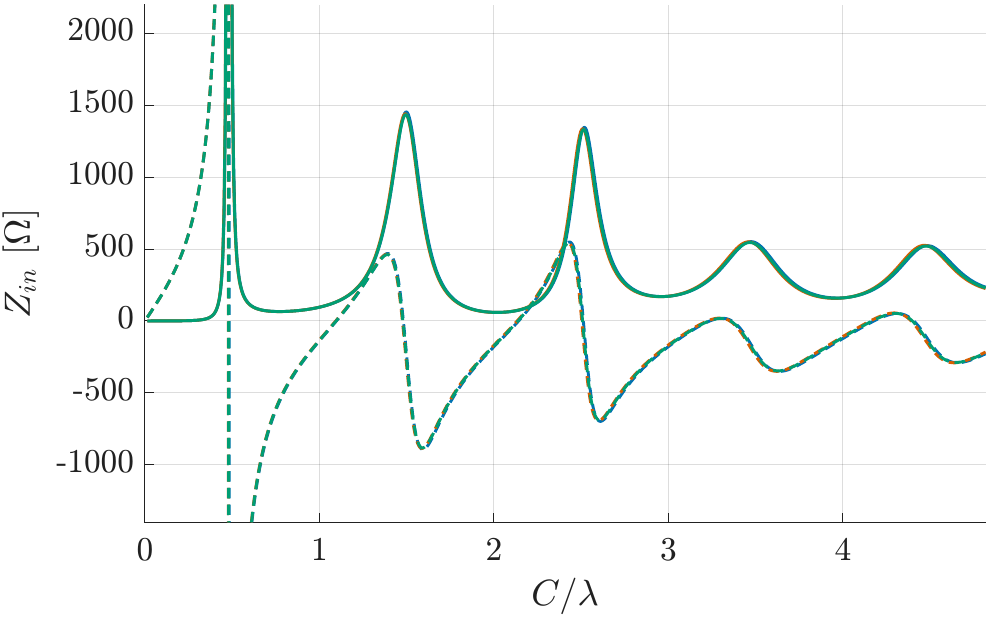}%
  \label{fig:square_Zin_BS2BS3}}\resultfigrowsep
\subfloat[$\delta_h^{(4)}$]{%
  \includegraphics[width=\resultfigwidth]{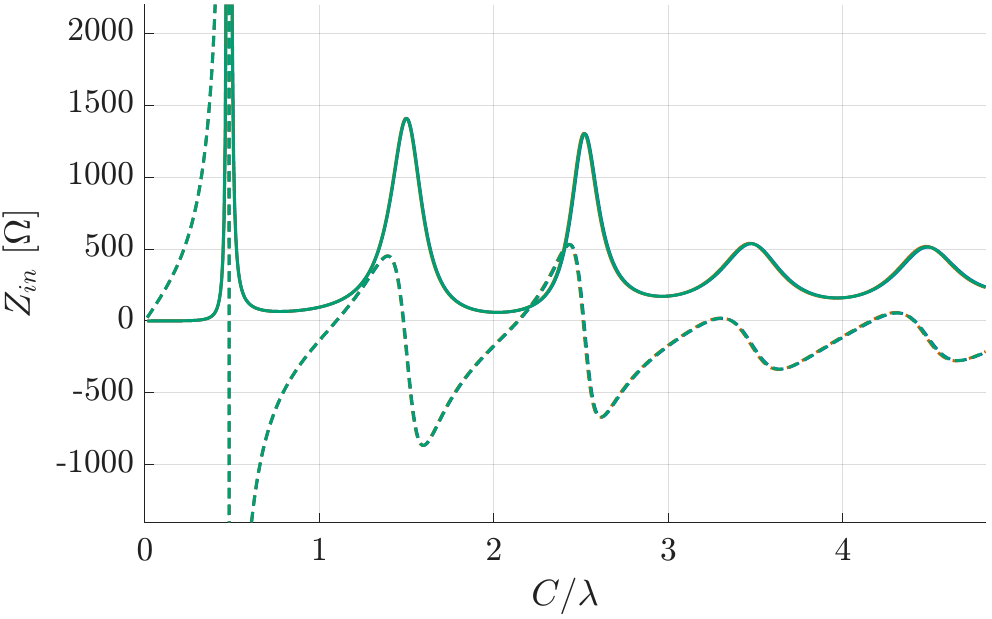}%
  \label{fig:square_Zin_BS4BS5}}\resultfigcolsep
\subfloat[$\delta_h^{\mathrm{iso}}$]{%
  \includegraphics[width=\resultfigwidth]{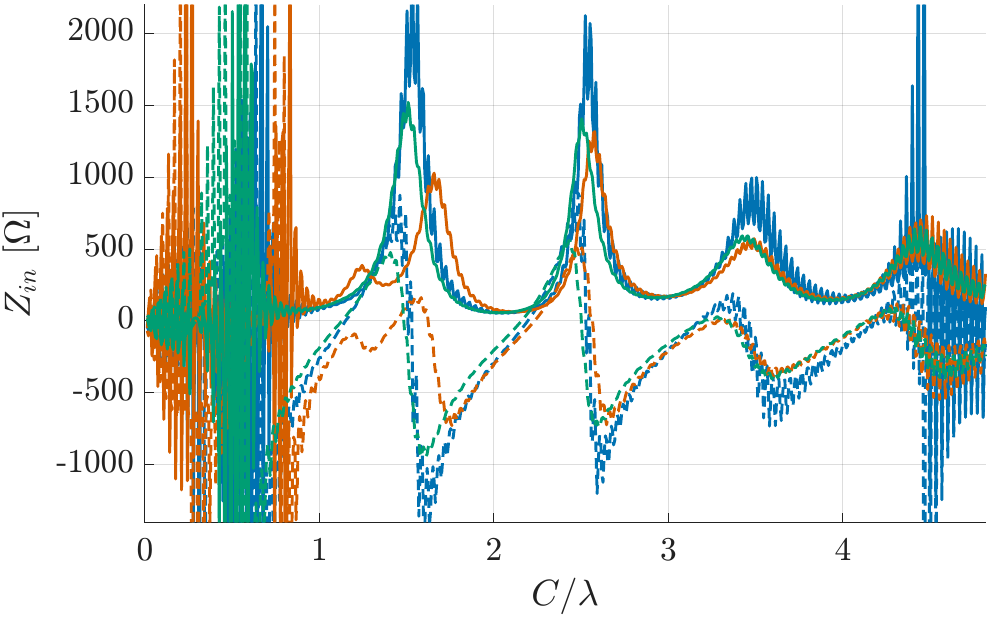}%
  \label{fig:square_Zin_iso}}
\caption{Input impedance of the square loop versus perimeter in wavelengths $C/\lambda$, with resistance $R_{\mathrm{in}}$ as solid lines and reactance $X_{\mathrm{in}}$ as dashed lines. Subpanels (a)--(d) use the four regularized delta functions $\delta_h^{(0)}$, $\delta_h^{(2)}$, $\delta_h^{(4)}$, and the isotropic $\delta_h^{\mathrm{iso}}$. Within each subpanel, three orientation curves are overlaid: axis-aligned $\hat{n}=\hat{z}$, face-diagonal $\hat{n}=(\hat{x}+\hat{y})/\sqrt{2}$, and body-diagonal $\hat{n}=(\hat{x}+\hat{y}+\hat{z})/\sqrt{3}$.}
\label{fig:square_Zin}
\end{figure}

\begin{figure}[t]
\centering
\subfloat[$\delta_h^{(0)}$]{%
  \includegraphics[width=\resultfigwidth]{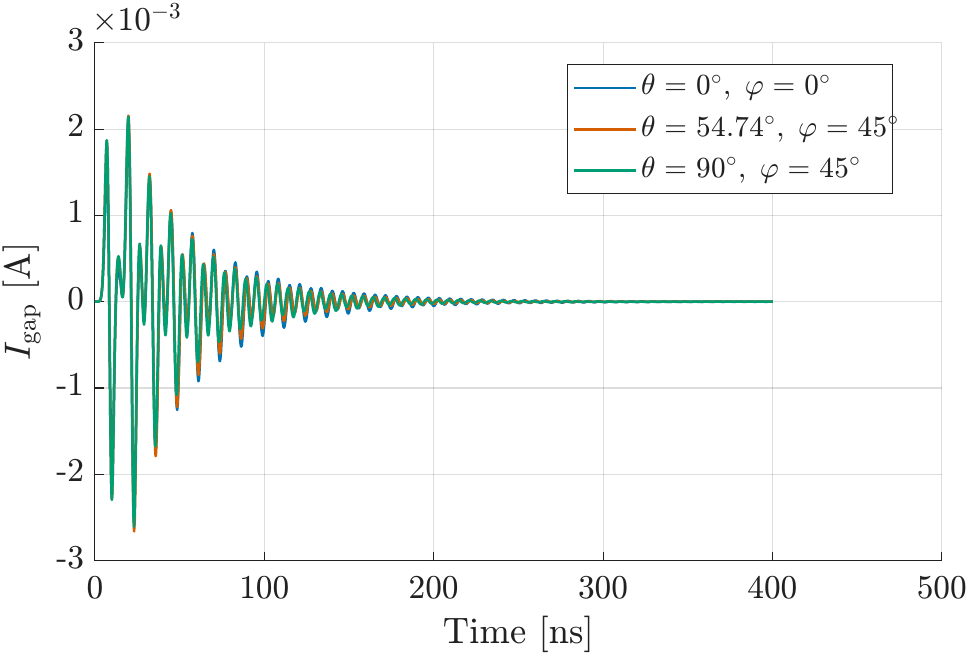}%
  \label{fig:square_Igap_t_BS0BS1}}\resultfigcolsep
\subfloat[$\delta_h^{(2)}$]{%
  \includegraphics[width=\resultfigwidth]{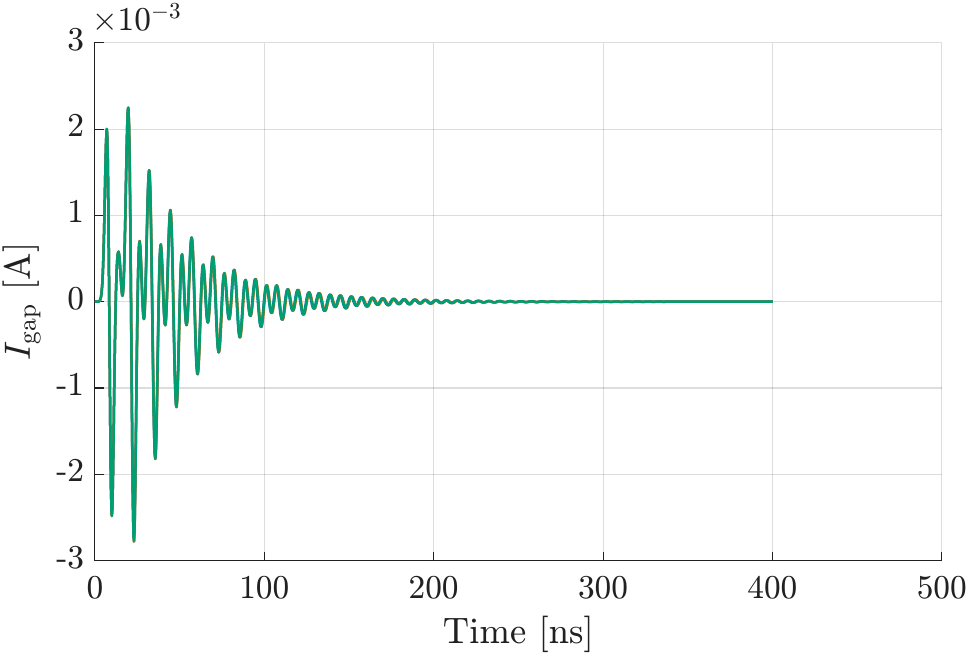}%
  \label{fig:square_Igap_t_BS2BS3}}\resultfigrowsep
\subfloat[$\delta_h^{(4)}$]{%
  \includegraphics[width=\resultfigwidth]{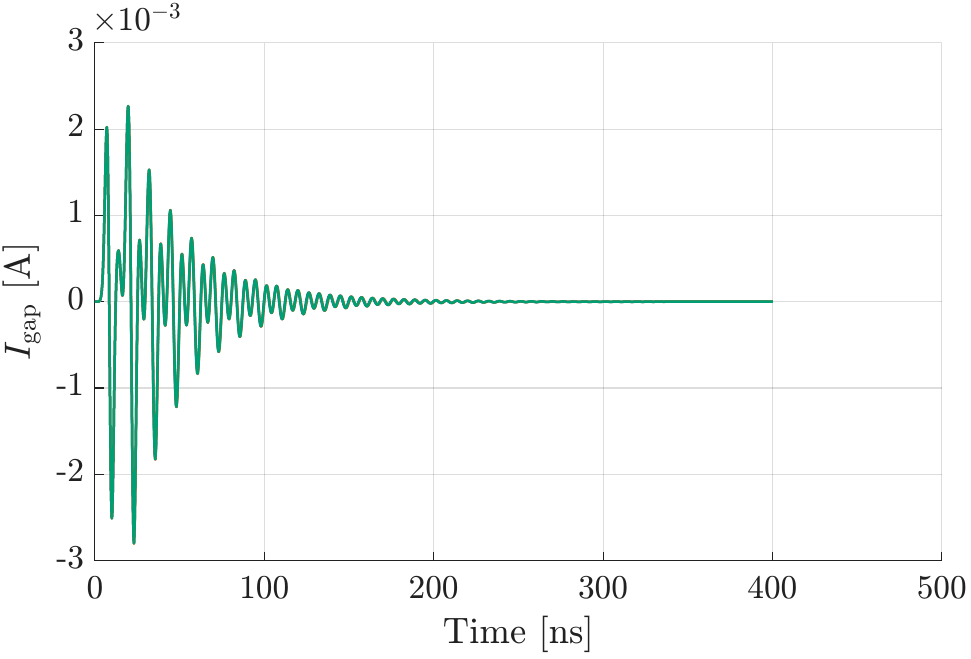}%
  \label{fig:square_Igap_t_BS4BS5}}\resultfigcolsep
\subfloat[$\delta_h^{\mathrm{iso}}$]{%
  \includegraphics[width=\resultfigwidth]{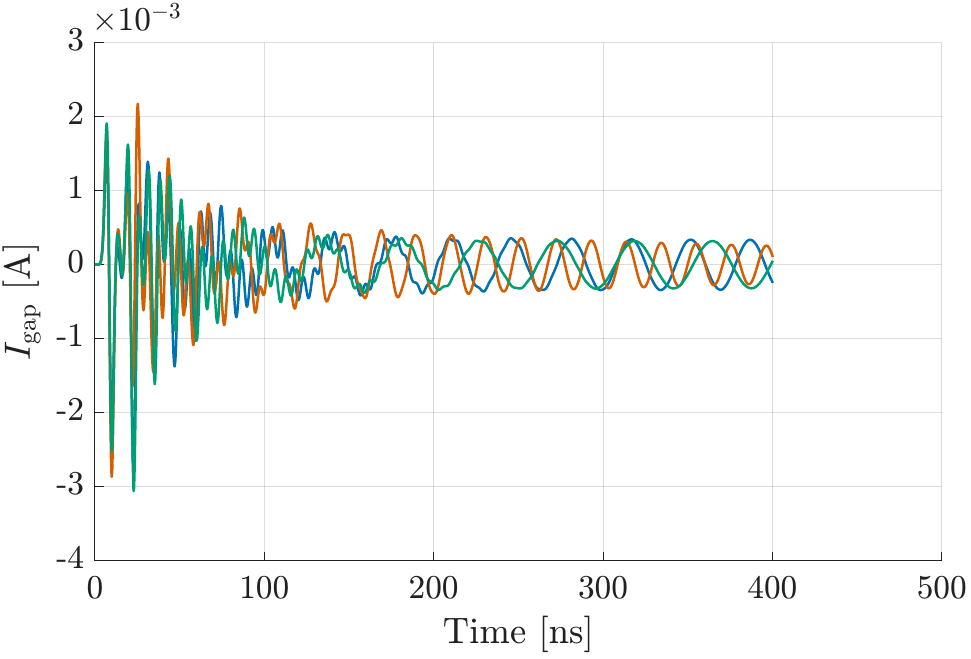}%
  \label{fig:square_Igap_t_iso}}
\caption{Gap-current time history $I_{\mathrm{gap}}(t)$ at the feed of the square loop. Subpanels (a)--(d) use the four regularized delta functions $\delta_h^{(0)}$, $\delta_h^{(2)}$, $\delta_h^{(4)}$, and the isotropic $\delta_h^{\mathrm{iso}}$, with the three orientation curves overlaid in each.}
\label{fig:square_Igap_t}
\end{figure}

Fig.~\ref{fig:square_Zin} reports the input impedance and Fig.~\ref{fig:square_Igap_t} the gap-current time histories, for the three orientations of the square loop antenna and four choices of regularized delta function.

The results for the square loop antenna are essentially analogous to those of the circular loop. The smoothest composite regularized delta functions tested, $\delta_h^{(2)}$ and $\delta_h^{(4)}$, produce orientation-independent impedance curves, and the corresponding gap currents ring down cleanly. The least regular composite regularized delta function, $\delta_h^{(0)}$, shows similar behavior, but with some mild variations across different orientations. The isotropic $\delta_h^{\mathrm{iso}}$ once again produces a persistent low-frequency current that contaminates the impedance with noise and renders it unphysical.

\section{Conclusion}
\label{sec:conclusion}

A closed wire carrying a constant current produces a current density whose divergence vanishes in the distributional sense. Charge conservation in thin-wire FDTD therefore reduces to a single requirement: the regularized current deposited on the Yee grid must inherit this property, so that the discrete divergence vanishes if the wire current is constant. This is the mathematical content of B\'{e}renger's geometric Q-cell criterion.

We satisfy this requirement with a family $\delta_h^{(n)}$ of composite B-spline regularized delta functions, using $\mathrm{BS}_n$ in the Yee-staggered direction and $\mathrm{BS}_{n+1}$ in the two node-aligned directions for each current component. The key algebraic property of cardinal B-splines is that a finite difference in the staggered direction can be expressed as a continuous derivative of the next-order B-spline, which allows the discrete divergence of a deposited constant current to be written as a total derivative in arc length that integrates to zero around any closed wire. This identity holds to machine precision as long as the panel integrals are evaluated exactly. Such exact evaluation is possible because the B-spline kernels are piecewise polynomial. In particular, along each piecewise-linear panel the integrand is a polynomial in arc length, and composite Gauss--Legendre quadrature with breakpoints at the wire's grid-plane crossings therefore integrates the panel integrals to machine precision. Without this piecewise-polynomial structure in either the regularized delta function or the discretization of the wire, an exact quadrature is not as readily available.

The interpolation operator is taken to be the discrete adjoint of the current deposition operator. By orthogonality of the gradient and curl subspaces in the discrete Helmholtz decomposition, an irrotational electric field then drives no net electromotive force around a closed wire loop, the discrete counterpart of the continuum identity that gradient fields have zero circulation. The coupled FDTD--wire system therefore neither accumulates spurious charge on the grid nor drives a spurious EMF around closed loops at the discrete level.

Numerical experiments on a center-fed dipole, a circular loop, and a square loop demonstrate that the composite regularized delta functions $\delta_h^{(0)}$, $\delta_h^{(2)}$, and $\delta_h^{(4)}$ produce input impedances consistent with the canonical responses of these antennas across all three orientations tested. The smoother members $\delta_h^{(2)}$ and $\delta_h^{(4)}$ collapse the orientation curves to within plotting tolerance, while $\delta_h^{(0)}$ shows a small residual orientation spread attributable to the lower regularity of $\mathrm{BS}_0$ in the staggered direction. The isotropic control $\delta_h^{\mathrm{iso}}$, which fails the charge-conservation condition, exhibits the persistent low-frequency current and orientation-dependent impedance discussed in Section~\ref{sec:results}.

A natural direction for future work is the extension of this framework from one-dimensional wires to two-dimensional surface antennas, such as microstrip patches and conducting sheets, for which the analogous compatibility condition requires the discrete divergence of a regularized surface current density to vanish on the Yee grid.

\section*{Acknowledgment}
Cole Gruninger is grateful for support from the Department of Defense (DoD) through the National Defense Science and Engineering Graduate (NDSEG) Fellowship Program. Boyce E. Griffith gratefully acknowledges support from the National Science Foundation (NSF) under grant number OAC 1931516.

%% file: supplementary_material_body.tex
We prove that, in the absence of absorbing boundary conditions, the leapfrog time-stepping scheme combined with the completely skew-symmetric spatial discretization of the thin-wire equations is energy preserving. For simplicity, we assume periodic boundary conditions on a finite Yee lattice, so that the $\ell^2$ inner products below are well-defined and finite. The argument adapts to free space provided the fields decay sufficiently rapidly at infinity that the corresponding sums converge absolutely.

We further assume, without loss of generality, that the spatial domain is a periodic unit cube discretized using $N$ grid points in each direction, with equal spacings $\dx = \dy = \dz = h = 1/N$. This is purely a notational simplification, and the argument extends to anisotropic grids with only minor changes. Before we begin the proof, we state some relevant definitions.

We work with two discrete $\ell^2$ inner products that the main text introduces for the Yee grid and the discretized wire. For grid-valued vector fields $\vF$ and $\vG$ on the staggered Yee lattice, the grid $\ell^2$ inner product is
\begin{align}
\langle \vF, \vG \rangle_h \equiv h^3 \sum_{i,j,k}\Big(&
F_x^{i+\tfrac{1}{2},j,k}\, G_x^{i+\tfrac{1}{2},j,k} \eqbreak
\eqcontinue{}+ F_y^{i,j+\tfrac{1}{2},k}\, G_y^{i,j+\tfrac{1}{2},k} \eqbreak
\eqcontinue{}+ F_z^{i,j,k+\tfrac{1}{2}}\, G_z^{i,j,k+\tfrac{1}{2}}\Big),
\label{eq:grid_inner_h}
\end{align}
in which $(i,j,k)$ ranges over $\{0,\dots,N-1\}^3$ with periodic identification at the boundaries. The subscript denotes the component of the vector field and is not meant to be a partial derivative. The same product applies to magnetic-field components $\vH$ evaluated at their corresponding staggered face locations $(i,j+\tfrac{1}{2},k+\tfrac{1}{2})$, $(i+\tfrac{1}{2},j,k+\tfrac{1}{2})$, $(i+\tfrac{1}{2},j+\tfrac{1}{2},k)$.

For wire-valued scalar quantities $A_q$, $B_q$ defined on the $N_w$ wire panels, the wire $\ell^2$ inner product is
\begin{equation}
\langle A, B \rangle_w \equiv \sum_{q=1}^{N_w} A_q\, B_q\, \Delta X_q,
\label{eq:wire_inner_sup}
\end{equation}
in which $\Delta X_q$ is the length of the $q$-th panel.

We next specify the discrete energy that is conserved by the scheme. Because the Yee scheme staggers solution values in time, we use an energy that combines two different time levels,
\begin{align}
\energy^n
={}& \frac{\epsilon_0}{2}\langle\vE^n,\vE^n\rangle_h \eqbreak
   \eqcontinue{}+ \frac{\mu_0}{2}\langle\vH^{n+1/2},\vH^{n-1/2}\rangle_h \eqbreak
   \eqcontinue{}+ \frac{C}{2}\langle V^n,V^n\rangle_w \eqbreak
   \eqcontinue{}+ \frac{L}{2}\langle I^{n+1/2},I^{n-1/2}\rangle_w.
\label{eq:energy_def}
\end{align}
To verify that $\energy^n$ is a valid energy, and hence non-negative, we apply the polarization identity to each time-staggered inner product to obtain
\begin{align}
\langle \vH^{n+1/2}, \vH^{n-1/2} \rangle_h
={}& \tfrac{1}{2}\|\vH^{n+1/2}\|_h^2 + \tfrac{1}{2}\|\vH^{n-1/2}\|_h^2 \eqbreak
   \eqcontinue{}- \tfrac{1}{2}\|\vH^{n+1/2} - \vH^{n-1/2}\|_h^2,
\end{align}
in which $\|\cdot\|_h$ is the norm induced by the discrete $\ell^2$ inner product on the grid, and likewise
\begin{align}
\langle I^{n+1/2}, I^{n-1/2} \rangle_w
={}& \tfrac{1}{2}\|I^{n+1/2}\|_w^2 + \tfrac{1}{2}\|I^{n-1/2}\|_w^2 \eqbreak
   \eqcontinue{}- \tfrac{1}{2}\|I^{n+1/2} - I^{n-1/2}\|_w^2,
\end{align}
in which $\|\cdot\|_w$ is the corresponding wire norm. The Yee scheme updates,
\begin{equation}
\vH^{n+1/2} = \vH^{n-1/2} -\frac{\dt}{\mu_0}\curlh \vE^n,
\end{equation}
and
\begin{equation}
I_q^{n+1/2} = I_q^{n-1/2} - \frac{\dt}{L}D_q V^n + \frac{\dt}{L}\Ical(\vE^n)_q,
\end{equation}
in which $D_q$ is the finite-difference operator $D_q V^n = (V_{q+1}^n - V_q^n)/\Delta X_q$, express the squared time differences in terms of $\vE^n$ and $V^n$. Substituting these expressions, the energy can be rewritten as
\begin{align}
\energy^n
={}& \frac{\epsilon_0}{2}\|\vE^n\|_h^2
   - \frac{\dt^2}{4\mu_0}\|\curlh \vE^n\|_h^2 \eqbreak
   \eqcontinue{}+ \frac{\mu_0}{4}\|\vH^{n+1/2}\|_h^2
   + \frac{\mu_0}{4}\|\vH^{n-1/2}\|_h^2 \nonumber\\
   &+ \frac{C}{2}\|V^n\|_w^2
   + \frac{L}{4}\|I^{n+1/2}\|_w^2 \eqbreak
   \eqcontinue{}+ \frac{L}{4}\|I^{n-1/2}\|_w^2
   - \frac{\dt^2}{4L}\|D_q V^{n} - \Ical(\vE^n)\|_w^2.
\end{align}
Applying Young's inequality $\|A - B\|_w^2 \le 2\|A\|_w^2 + 2\|B\|_w^2$ to the last term yields
\begin{align}
\energy^n
\geq{}& \frac{\epsilon_0}{2}\|\vE^n\|_h^2
   - \frac{\dt^2}{4\mu_0}\|\curlh \vE^n\|_h^2 \eqbreak
   \eqcontinue{}+ \frac{\mu_0}{4}\|\vH^{n+1/2}\|_h^2
   + \frac{\mu_0}{4}\|\vH^{n-1/2}\|_h^2 \nonumber\\
   &+ \frac{C}{2}\|V^n\|_w^2
   + \frac{L}{4}\|I^{n+1/2}\|_w^2 \eqbreak
   \eqcontinue{}+ \frac{L}{4}\|I^{n-1/2}\|_w^2
   - \frac{\dt^2}{2L}\|D_q V^{n}\|_w^2 \eqbreak
   \eqcontinue{}- \frac{\dt^2}{2L}\|\Ical(\vE^n)\|_w^2.
\end{align}
Bounding $\|\curlh\vE^n\|_h$, $\|D_q V^n\|_w$, and $\|\Ical(\vE^n)\|_w$ by the corresponding induced operator norms then gives
\begin{align}
\energy^n
\geq{}& \frac{\epsilon_0}{2}\|\vE^n\|_h^2
   - \frac{\dt^2}{4\mu_0}\|\curlh\|^2_{\mathrm{op}}\|\vE^n\|_h^2 \eqbreak
   \eqcontinue{}+ \frac{\mu_0}{4}\|\vH^{n+1/2}\|_h^2
   + \frac{\mu_0}{4}\|\vH^{n-1/2}\|_h^2 \nonumber\\
   &+ \frac{C}{2}\|V^n\|_w^2
   + \frac{L}{4}\|I^{n+1/2}\|_w^2 \eqbreak
   \eqcontinue{}+ \frac{L}{4}\|I^{n-1/2}\|_w^2
   - \frac{\dt^2}{2L}\|D_q\|^2_{\mathrm{op}}\|V^{n}\|_w^2 \eqbreak
   \eqcontinue{}- \frac{\dt^2}{2L}\|\Ical\|^2_{\mathrm{op}}\|\vE^n\|_h^2,
\end{align}
in which $\|\cdot\|_{\mathrm{op}}$ denotes the operator norm induced by the relevant inner product. Collecting like terms, we have
\begin{align}
\energy^n
\geq{}& \left(\frac{\epsilon_0}{2}
   - \frac{\dt^2}{4\mu_0}\|\curlh\|^2_{\mathrm{op}}
   - \frac{\dt^2}{2L}\|\Ical\|^2_{\mathrm{op}}\right)\|\vE^n\|_h^2 \nonumber\\
   &+ \frac{\mu_0}{4}\|\vH^{n+1/2}\|_h^2
   + \frac{\mu_0}{4}\|\vH^{n-1/2}\|_h^2 \eqbreak
   \eqcontinue{}+ \left(\frac{C}{2}
   - \frac{\dt^2}{2L}\|D_q\|^2_{\mathrm{op}}\right)\|V^n\|_w^2 \nonumber\\
   &+ \frac{L}{4}\|I^{n+1/2}\|_w^2
   + \frac{L}{4}\|I^{n-1/2}\|_w^2.
\end{align}
Choosing the time-step size $\dt$ small enough that the terms in parentheses are non-negative therefore guarantees $\energy^n \ge 0$. The first term constrains the electromagnetic time step,
\begin{equation}
\dt^2 \le \frac{2\mu_0\epsilon_0}{\|\curlh\|^2_{\mathrm{op}} + (2\mu_0/L)\|\Ical\|^2_{\mathrm{op}}},
\label{eq:dt_em}
\end{equation}
and the second constrains the wire time step,
\begin{equation}
\dt^2 \le \frac{CL}{\|D_q\|^2_{\mathrm{op}}}.
\label{eq:dt_wire}
\end{equation}
In the absence of wire coupling, $\|\Ical\|_{\mathrm{op}}=0$, and \eqref{eq:dt_em} reduces to $\dt^2\|\curlh\|^2_{\mathrm{op}}\le 2\mu_0\epsilon_0$, which recovers the standard CFL constraint for the leapfrog FDTD scheme with wave speed $c=1/\sqrt{\mu_0\epsilon_0}$. The wire constraint \eqref{eq:dt_wire} is the analogous transmission-line CFL condition with $\dt=\mathcal{O}(\Delta X/c)$ and $c=1/\sqrt{LC}$.

The coupling of the wire to the background electromagnetic field is captured by $\|\Ical\|^2_{\mathrm{op}}$ and tightens the time-step size constraint relative to the pure-FDTD CFL condition. The exact value of $\|\Ical\|_{\mathrm{op}}$ depends on the choice of regularized delta function and on the placement and orientation of the wire relative to the background grid. Its scaling in $h$ can be anticipated using the Cauchy-Schwarz inequality combined with the zeroth-moment identity for the regularized delta function, $h^3\sum_{i,j,k}\delta_h=1$. Writing the interpolation operator as
\begin{equation}
\Ical(\vE)_q = h^3\sum_{i,j,k}\delta_h(\vx_{ijk}-\vX_q)\,\vE(\vx_{ijk})\cdot\vtau_q,
\end{equation}
the Cauchy-Schwarz inequality and the zeroth-moment identity yield the bound
\begin{equation}
|\Ical(\vE)_q|^2 \le h^3\sum_{i,j,k}\delta_h(\vx_{ijk}-\vX_q)|\vE(\vx_{ijk})|^2.
\end{equation}
Multiplying each side by $\Delta X_q$, summing in $q$, and exchanging the order of summation gives
\begin{align}
\|\Ical(\vE)\|_w^2
\le{}& h^3\sum_{i,j,k}|\vE(\vx_{ijk})|^2 \eqbreak
\eqcontinue{}\times\left(\sum_q \delta_h(\vx_{ijk}-\vX_q)\,\Delta X_q\right).
\end{align}
The sum over the wire degrees of freedom in parentheses is a discrete one-dimensional regularization of the wire's distributional source term. In the regime in which the panel and grid spacings are comparable, $\Delta X \sim h$, the regularized delta function $\delta_h \sim h^{-3}$ integrated along a one-dimensional wire contributes one factor of $h$, giving $\|\Ical\|_{\mathrm{op}} = \mathcal{O}(h^{-1})$, or equivalently $\|\Ical\|_{\mathrm{op}} = \mathcal{O}(\Delta X^{-1})$. This reflects the singular nature of restricting a three-dimensional field to a one-dimensional curve, which is unbounded at the continuum level. For a codimension-one surface, integrating $\delta_h$ over a two-dimensional support contributes two factors of $h$, yielding the weaker scaling $\|\Ical\|_{\mathrm{op}} = \mathcal{O}(\Delta X^{-1/2})$.

We remark that the interpolation operator used in the numerical experiments of the main text is implemented via composite Gauss--Legendre quadrature over each wire panel rather than the point evaluation written above. Because the panel quadrature preserves the zeroth-moment identity and only refines the integration along the wire's tangent direction, the argument carries over with at most an $\mathcal{O}(1)$ change in the constant, and the scaling $\|\Ical\|_{\mathrm{op}}=\mathcal{O}(\Delta X^{-1})$ should apply equally to the quadrature-based interpolation.

Having shown that the energy is non-negative under a suitable time-step size constraint, we next show that taking the wire interpolation and current-deposition operators to be discrete adjoints of one another ensures that this energy is exactly conserved.

To do this, we analyze the difference $\energy^{n+1}-\energy^n$ and show that it is identically zero if $\Scal$ and $\Ical$ are discrete adjoints. The leapfrog scheme provides the updates
\begin{align}
\vE^{n+1} &= \vE^n + \frac{\dt}{\epsilon_0}\left(\curlh \vH^{n+1/2} - \vJ^{n+1/2}\right), \\
\vH^{n+3/2} &= \vH^{n+1/2} - \frac{\dt}{\mu_0}\curlh \vE^{n+1}, \\
V_q^{n+1} &= V_q^n - \frac{\dt}{C}D_q I^{n+1/2}, \\
I_q^{n+3/2} &= I_q^{n+1/2} - \frac{\dt}{L}D_q V^{n+1} \eqbreak
            \eqcontinuequad{}+ \frac{\dt}{L}\Ical(\vE^{n+1})_q,
\end{align}
in which $\vJ^{n+1/2} = \Scal(I_q^{n+1/2})$ is the source current density deposited onto the grid by the wire current deposition operator. The staggered placement of $V$ and $I$ on the wire lets the same operator $D_q$ act on both, with the skew-adjoint identity
\begin{equation}
\langle D_q V, I\rangle_w = -\langle V, D_q I\rangle_w,
\label{eq:wire_skew_adjoint}
\end{equation}
which follows from summing by parts and using the periodic boundary conditions.

For simplicity, we treat the electromagnetic and wire contributions to the energy difference separately. Beginning with the electromagnetic piece, we have
\begin{align}
\eqfirstlhs{\frac{\epsilon_0}{2}\left(\|\vE^{n+1}\|_h^2-\|\vE^n\|_h^2\right)} \eqfirststep
\eqfirstrel \frac{\epsilon_0}{2}\langle \vE^{n+1}+\vE^n,\vE^{n+1}-\vE^n\rangle_h \nonumber\\
\eqrel \frac{\dt}{2}\langle \vE^{n+1}+\vE^n,\curlh\vH^{n+1/2}\rangle_h \eqbreak
   \eqcontinue{}- \frac{\dt}{2}\langle \vE^{n+1}+\vE^n,\vJ^{n+1/2}\rangle_h.
\end{align}
The magnetic piece follows from subtracting the half-step updates at $n+1$ and $n$, which gives $\vH^{n+3/2}-\vH^{n-1/2}=-\tfrac{\dt}{\mu_0}\curlh(\vE^{n+1}+\vE^n)$, so that
\begin{align}
\eqfirstlhs{\frac{\mu_0}{2}\left(\langle\vH^{n+3/2},\vH^{n+1/2}\rangle_h-\langle\vH^{n+1/2},\vH^{n-1/2}\rangle_h\right)} \eqfirststep
\eqfirstrel \frac{\mu_0}{2}\langle\vH^{n+1/2},\vH^{n+3/2}-\vH^{n-1/2}\rangle_h \nonumber\\
\eqrel -\frac{\dt}{2}\langle\vH^{n+1/2},\curlh(\vE^{n+1}+\vE^n)\rangle_h \nonumber\\
\eqrel -\frac{\dt}{2}\langle\curlh\vH^{n+1/2},\vE^{n+1}+\vE^n\rangle_h,
\end{align}
in which the last equality uses summation by parts. Adding the two contributions, the curl terms cancel and the electromagnetic contribution to $\energy^{n+1}-\energy^n$ reduces to
\begin{equation}
-\frac{\dt}{2}\langle \vE^{n+1}+\vE^n,\vJ^{n+1/2}\rangle_h.
\label{eq:em_delta}
\end{equation}
The contribution due to the wire is handled analogously. Using $V^{n+1}-V^n=-\tfrac{\dt}{C}D_q I^{n+1/2}$ together with the skew-adjoint identity \eqref{eq:wire_skew_adjoint},
\begin{align}
\eqfirstlhs{\frac{C}{2}\left(\|V^{n+1}\|_w^2-\|V^n\|_w^2\right)} \eqfirststep
\eqfirstrel \frac{C}{2}\langle V^{n+1}+V^n, V^{n+1}-V^n\rangle_w \nonumber\\
\eqrel -\frac{\dt}{2}\langle V^{n+1}+V^n, D_q I^{n+1/2}\rangle_w \nonumber\\
\eqrel \frac{\dt}{2}\langle D_q(V^{n+1}+V^n), I^{n+1/2}\rangle_w,
\end{align}
and the corresponding calculation for the currents, with $I^{n+3/2}-I^{n-1/2}=-\tfrac{\dt}{L}D_q(V^{n+1}+V^n)+\tfrac{\dt}{L}\Ical(\vE^{n+1}+\vE^n)$, yields
\begin{align}
\eqfirstlhs{\frac{L}{2}\left(\langle I^{n+3/2},I^{n+1/2}\rangle_w-\langle I^{n+1/2},I^{n-1/2}\rangle_w\right)} \eqfirststep
\eqfirstrel \frac{L}{2}\langle I^{n+1/2},I^{n+3/2}-I^{n-1/2}\rangle_w \nonumber\\
\eqrel -\frac{\dt}{2}\langle I^{n+1/2},D_q(V^{n+1}+V^n)\rangle_w \eqforcedbreak
   \eqforcedcontinue{}+ \frac{\dt}{2}\langle I^{n+1/2},\Ical(\vE^{n+1}+\vE^n)\rangle_w.
\end{align}
Adding these, the $D_q$ terms cancel, and the wire contribution to $\energy^{n+1}-\energy^n$ reduces to
\begin{equation}
\frac{\dt}{2}\langle I^{n+1/2},\Ical(\vE^{n+1}+\vE^n)\rangle_w.
\label{eq:wire_delta}
\end{equation}
Combining \eqref{eq:em_delta} and \eqref{eq:wire_delta}, we obtain
\begin{align}
\energy^{n+1}-\energy^n
={}& -\frac{\dt}{2}\langle \vE^{n+1}+\vE^n,\vJ^{n+1/2}\rangle_h \eqbreak
   \eqcontinue{}+ \frac{\dt}{2}\langle I^{n+1/2},\Ical(\vE^{n+1}+\vE^n)\rangle_w.
\end{align}
Substituting $\vJ^{n+1/2}=\Scal(I_q^{n+1/2})$ shows that both terms pair the wire current $I^{n+1/2}$ with the field combination $\vE^{n+1}+\vE^n$, one through the deposition operator and the other through the interpolation operator. If $\Scal$ and $\Ical$ are chosen to be discrete adjoints in the sense that
\begin{equation}
\langle \vF,\Scal(J_q)\rangle_h = \langle J,\Ical(\vF)\rangle_w
\label{eq:adjoint_sup}
\end{equation}
for all grid vector fields $\vF$ and wire currents $J$, then the two terms agree and
\begin{equation}
\energy^{n+1} = \energy^n.
\end{equation}
The leapfrog FDTD--thin-wire system therefore conserves the discrete energy $\energy^n$ exactly under periodic boundary conditions, provided the wire-to-grid current deposition and grid-to-wire interpolation operators are chosen as discrete adjoints with respect to the grid and wire $\ell^2$ inner products.